\documentclass[11 pt]{amsart}
\usepackage{geometry}                
\usepackage{amsmath, amsthm, amssymb} 
\usepackage{graphicx}
\usepackage{amssymb}
\usepackage{epstopdf}
\usepackage{float}
\newtheorem{corollary}{Corollary}
\theoremstyle{definition}
\newtheorem{definition}{Definition}
\theoremstyle{remark}
\newtheorem{remark}{Remark}

\title[Simplicial Ricci Flow]{Simplicial Ricci Flow}

\author[Miller, McDonald Alsing, Gu \& Yau]{Warner A. Miller,${}^{1,2}$ Jonathan R. McDonald,${}^1$ Paul M. Alsing,${}^3$  \\David Gu$^{4}$ 
 \ \& \ Shing-Tung Yau${}^1$}
\address{\dag Department of Mathematics, Harvard University, Cambridge MA 02138.}
\address{\ddag Department of Physics, Florida Atlantic University, Boca Raton, FL 33431} 
\address{${}^*$ Air Force Research Laboratory, Information Directorate, Rome, NY 13441 } 
\address{${}^{**}$ Department of Computer Science, Stony Brook University, Stony Brook, NY 11794 }

\date{February 4, 2013}                                           

\begin{document}

\begin{abstract}
We construct a discrete form of Hamilton's Ricci flow (RF) equations for a $d$-dimensional piecewise flat simplicial geometry, ${\mathcal S}$.  These new algebraic equations are derived using the discrete formulation of Einstein's theory of general relativity known as Regge calculus.   A  Regge-Ricci flow (RRF) equation can be associated to each edge, $\ell$,  of a  simplicial lattice.  In defining this equation, we find it convenient to utilize both the simplicial lattice ${\mathcal S}$ and its circumcentric dual lattice, ${\mathcal S}^*$. In particular, the RRF equation associated to $\ell$ is naturally defined on a $d$-dimensional hybrid block connecting $\ell$ with its $(d-1)$-dimensional circumcentric dual cell, $\ell^*$.  We show that this equation  is expressed as the proportionality between (1) the simplicial Ricci tensor, $Rc_\ell$,  associated with the edge $\ell\in{\mathcal S}$, and (2) a certain volume weighted average of the fractional rate of change of the edges, $\lambda\in \ell^*$, of the circumcentric dual lattice, ${\mathcal S}^*$,  that are in the dual of $\ell$. The inherent orthogonality between elements of $\mathcal S$ and their duals in ${\mathcal S}^*$ provide a simple geometric representation of Hamilton's RF equations.  In this paper we utilize the well established theories of Regge calculus, or equivalently discrete exterior calculus, to construct these equations. We solve these equations for a few illustrative examples. 
 \ \\ \\ {E-mail: wam@fau.edu}
\end{abstract}

\maketitle

\section{Ricci Flow on Simplicial Geometries}
\label{sec:rfsg}
Hamilton's Ricci flow (RF)  technique has profoundly impacted the mathematical sciences and engineering fields \cite{Cao:2003,Chow:2004,Chow:2006,Chow:2007}.  While RF has been used to solve the Poincar\'e conjecture \cite{Perelman:2003},  the bulk of the applications of this technique have been limited to 2--dimensional geometries \cite{Gu:2012}.  It is well established that the analysis of a geometry with complex topology is most naturally represented as a coordinate--free, piecewise--flat simplicial lattice, e.g. as is routinely done using finite volume \cite{Peiro:2005}, finite element \cite{Humphries:1997}, Regge calculus \cite{Regge:1961,Gentle:1998} and discrete exterior calculus  \cite{Hirani:2005}.  In coordinate--based approaches, the difficulties in properly handling coordinate singularities, as well as the joining of multiple coordinate charts over the manifold, can be prohibitive.  On the other hand,  the utility of piecewise--flat simplicial geometries in analyzing the RF of 2--dimensional geometries is well established and proven to be effective \cite{Gu:2012,Chow:2003}.   This utility should extend to three and higher dimensions, and preliminary work is underway in this direction \cite{Glickenstein:2011,G:2005,Ge:2013,Forman:2003}. 

Simplicial and combinatorial RF is well established in 2--dimensions; however, we are unaware of any formulation of the RF equations for a piecewise-flat simplicial geometry in dimensions  greater than or equal to 3.  It is the purpose of this manuscript to develop such higher--dimensional Regge-Ricci flow (RRF) equations.  These  equations can be used in the analysis of topology and geometry, both (1) in the numerical analysis of RF and, perhaps more importantly, (2)  as a mathematical tool to bound the Ricci curvature for discrete geometries, and to analyze and handle higher--dimensional RF singularities \cite{LinYau:2010,Knopf:2009}.  The topological taxonomy afforded by RF is richer when transferring from 2 to 3--dimensions. In particular,  the uniformization theorem says that any  2--geometry will evolve under RF to a constant curvature sphere, plane or hyperboloid, while in 3--dimensions any manifold will diffuse into a connected sum of prime manifolds \cite{Thurston:1997}. 

We consider in this paper the structure of Hamilton's RF equations for a piecewise-flat simplicial geometry, ${\mathcal S}$, and  its circumcentric dual lattice, ${\mathcal S}^*$.  In Sec.~\ref{sec:theorem} we use the simplicial geometry approach known as  Regge calculus (RC) or its equivalent, discrete exterior calculus (DEC), to motivate our definition of a discrete representation of Hamilton's Ricci flow equation on each edge, $\lambda$, of the circumcentric dual lattice, ${\mathcal S}^*$.  We will refer to these as the dual-edge RRF equations.  Based on this definition, and using the results of   \cite{McDonald:2012,M:1997,MM:2008} on the analysis of curvature on piecewise-flat geometries,  we prove a corollary re-expressing these dual-edge RRF equations as a corresponding equation for each edge, $\ell$, of the $d$--dimensional piecewise flat simplicial geometry, $\mathcal S$.  We refer to these equations as the simplicial RRF equations, or simply as the RRF equations. 
 
\begin{definition}
\label{T1}
We define the dual-edge Regge-Ricci flow equation for any compact, piecewise--flat simplicial geometry, $\mathcal S$, as an equation for each edge, $\lambda$, in the circumcentric dual lattice, ${\mathcal S}^*$,
\begin{equation}
\left( \begin{array}{c}Fractional\ rate\ of\ change\\ of\  a\  circumcentric\\ dual\ edge,\\  \frac{1}{\lambda} \frac{\partial \lambda}{\partial t}\\   \end{array}  \right) = 
-\left( \begin{array}{c}Ricci\ tensor\ associated\\ to\ the\ circumcentric\\ dual\ edge,\\ Rc_\lambda    \end{array}  \right).
\end{equation}
\end{definition}

\begin{corollary}
\label{C1}
The simplicial Regge-Ricci flow equation for any compact, piecewise--flat simplicial geometry, $\mathcal S$, can be expressed as an equation for each edge, $\ell \in {\mathcal S}$,
\begin{equation}
\label{Eq:srfell}
\left( \begin{array}{c}The\  volume\ averaged\\  fractional\ rate\ of\ change\\ of\  the\ dual\ edge,\\  \left< \frac{1}{\lambda} \frac{\partial \lambda}{\partial t}\right>_{\ell}   \end{array}  \right) = 
-\left( \begin{array}{c}Ricci\ tensor\ associated\ to\\ the\ simplicial\ edge,\\ Rc_\ell    \end{array}  \right).
\end{equation} 
\end{corollary}

We express Hamilton's RF equation in an orthogonal vielbein basis in Sec.~\ref{sec:rf} so as to introduce the proper background and motivation for the definition of the RRF equations.  In Sec.~\ref{sec:curvature} we provide an overview of the structure of the Riemann tensor ($Rm$), Ricci tensor ($Rc$) and scalar curvature ($R$) on a discrete piecewise--flat simplicial geometry.  In Sec.~\ref{sec:theorem} we motivate Def.~1 and prove Cor.~1.  In Sec.~\ref{sec:RRF3}, we explicitly construct the 3--dimensional RRF equations and solve these analytically for two well known analytic examples, the 3--sphere and the 3--cylinder geometries. We conclude this manuscript in Sec.~\ref{Sec:Fini} by addressing a few important open problems.

\section{Hamilton's Ricci Flow Equation}
\label{sec:rf}
The RF equation of a $d$-dimensional Riemannian geometry was first introduced by Hamilton \cite{Hamilton:1982},   
\begin{equation}
\label{Eq:RF}
\frac{\partial g_{\mu\nu}}{\partial t} = -2\, R_{\mu\nu},
\end{equation}
where $R_{\mu\nu}$ are the covariant components of the Ricci tensor $Rc$, and $g_{\mu\nu}$ are the components of the metric tensor.
It is our goal to reconstruct the  RF equation of Hamilton on a piecewise--flat simplicial geometry, ${\mathcal S}$, by using the formalism of RC 
\cite{Regge:1961}.  To this end, and without loss of generality, consider a point $p \in \left\{{\mathcal M},g\right\}$ where we can construct a time--independent orthogonal basis, $\left\{ {\pmb e}_a \right\},\  \forall a = 1, 2 \ldots. d$,  i.e.
\begin{equation}
\label{Eq:edot}
\dot{\pmb e}_a \equiv \frac{\partial\,{\pmb e}_a}{\partial t} = 0.
\end{equation}
The contravariant components of the metric tensor in this basis,
\begin{equation}
{\pmb g} = {\pmb e_a} g^{ab} {\pmb e_b},
\end{equation}
will be, by construction, diagonal,
\begin{equation}
g^{ab} = diag\left\{ g^{aa} \right\}.
\end{equation}
In this specially-chosen basis (Eq.~\ref{Eq:edot}) the time rate of change of the metric will be expressed solely in terms of the time derivatives of its $d$ diagonal components. 
\begin{equation}
\dot {\pmb g} = {\pmb e_a}\, \dot{g}^{ab} \,{\pmb e_b} = {\pmb e_a}\, \dot{g}^{aa}\, {\pmb e_a}
\end{equation}
We will find it convenient for our purposes  to express the RF equation (Eq.~\ref{Eq:RF})  in its  mixed--index  form, i.e. as the components of a vector-valued 1-form equation in this basis at point $p$ in the sense of Cartan geometry \cite{Cartan:1951}.    The components of these mixed--index RF equations take a simple form, 
\begin{equation}
\label{eq:crf}
g^{aa}\, \dot g_{aa} = -2\, R^{a}{}_{a}.
\end{equation}
Since the metric at point $p$ is diagonal then, $g^{aa} = 1/g_{aa}$, and the RF equation takes a simple form, 
\begin{equation}
\label{Eq:frcc}
\left( \begin{array}{c} Fractional\ rate\ of\\ change\ of\ g_{aa}\end{array}\right)=\frac{{\dot g}_{aa}}{g_{aa}} = -2 R^a{}_a.
\end{equation}

\section{Curvature of a Piecewise--Flat $d$-Dimensional Simplicial Geometry}
\label{sec:curvature}

In this section we provide a brief review of the curvature on a piecewise--flat simplicial lattice with emphasis on its application  to the RRF equation.  Although a more complete and thorough description will appear in the literature \cite{McDonald:2012}, we felt that a review would make the derivation more complete and provide a slightly different emphasis needed for the derivation. In particular, we will outline the construction of the  $Rm$,  $Rc$ and $R$ on the simplicial lattice as well as  its circumcentric dual lattice.  We find it convenient to leverage the framework of  RC \cite{Regge:1961,FL:1984,M:1986,M:1997,MM:2008,AMM:2011}.  However, an essentially equivalent treatment of curvature on piecewise--flat simplicial manifolds is  afforded  by recent developments in DEC \cite{Hirani:2005}.  Both approaches, RC and DEC, appear to have independently encoded the geometry and curvature structure on discrete simplicial lattices in similar ways, and both utilize the circumcentric dual lattice in a fundamentally important way \cite{McDonald:2012}.  

In the Euclidean form of RC the geometry of the interior of  each of its $ d $--dimensional simplicies is assumed to be flat Euclidean space, $E^d$. The geometry of each simplex in $ \mathcal S $ is completely determined by the square of its ${d\choose2} = \frac{d(d-1)}{2}$ edge lengths.  Therefore, the entire  lattice geometry and its curvature structure is completely determined by the squared length of all its  edges, $\left\{ \ell^2_i \right\}$, and only these edges.   The set of all  edge lengths of the simplicial geometry is the simplicial analogue of the metric, i.e.  
\begin{equation}
{\pmb g} \longleftrightarrow \{\ell_i^2\},
\end{equation}
 in that any geometric quantity on the lattice, e.g. its lengths, areas, volumes, curvatures, and even the dual lattice are defined solely in terms of these edge lengths and nothing else.   
 
 RC was originally developed for general relativity as a piecewise--flat, coordinate--free simplicial representation of the 4-dimensional spacetime geometry.  The  geometry interior to each simplex  is assumed to be a  Minkowski flat spacetime.  However,  this theory is equally applicable for the curvature analysis of a  $d$--dimensional simplicial spacelike geometry of interest here.  Accordingly, we will assume that the geometry interior to each $d$--simplex is Euclidean flat space.  In this transition, the formalism is otherwise unaltered  from 4-dimensional spacetime geometries to purely spatial geometries.  We will use this fact to implement RC and to construct the RRF equation in $d$--dimensions in this manuscript.

\begin{figure}
\centering
\includegraphics[height=6.0in]{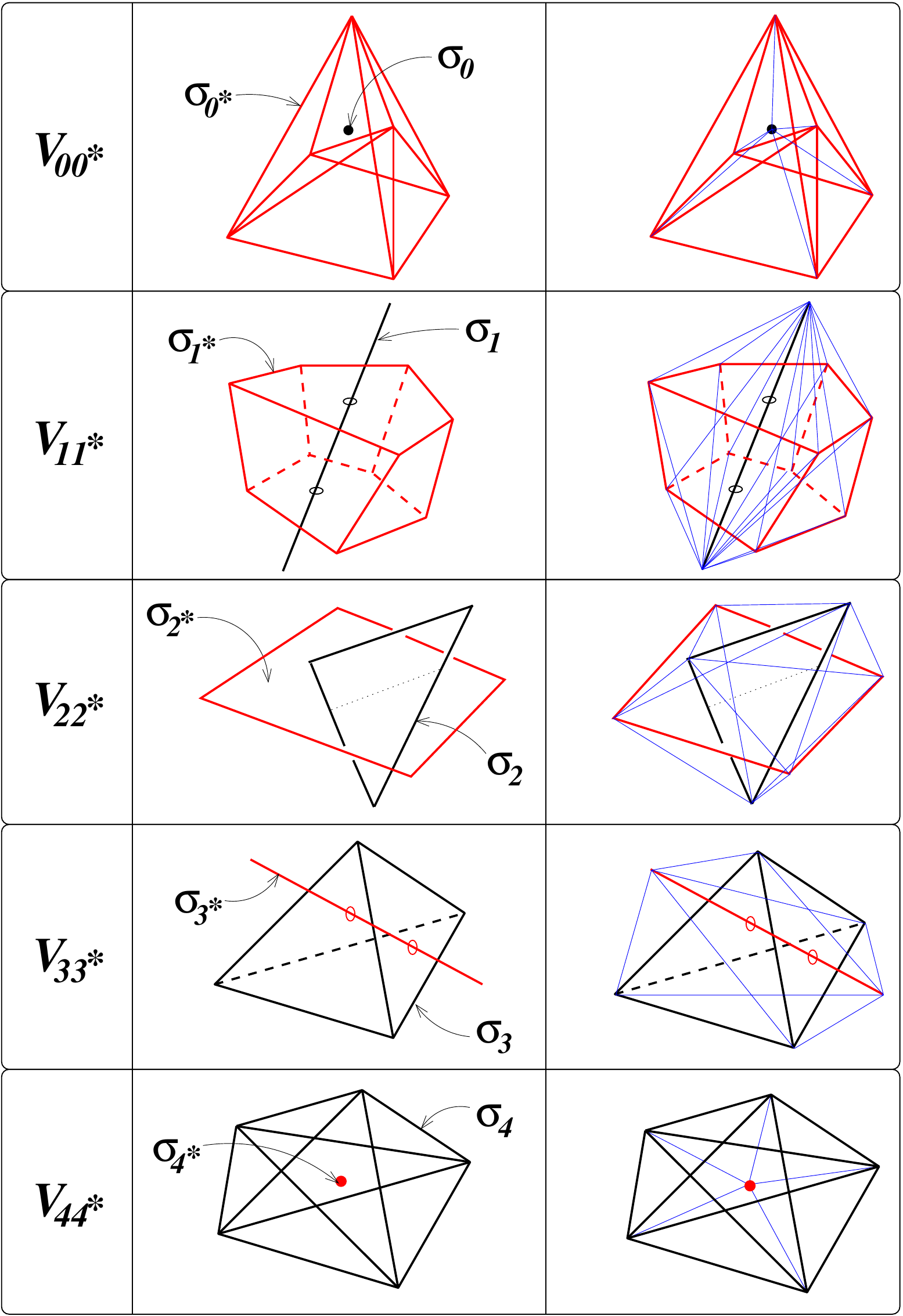}
\caption{ Each of the five rows in this table illustrates one of the  five hybrid 4--dimensional polytopes.   The black edges (thick lines)  in the diagram are in the simplicial Delaunay lattice, ${\mathcal S}$ while the red (thick grey)  edges are in the circumcentric dual Voronoi lattice, ${\mathcal S}^*$.  The blue (thin)  edges are added in the right--most column and connect the vertices of the Delaunay lattice to the vertices of the Voronoi lattice thus forming a 4--dimensional polytope.  Each of these 5 polytopes, when applied throughout all similar elements of $\mathcal S$, provide a complete tessellation of  the  lattice spacetime. The symbols $V_{kk^*}$ in the left column completely characterize the hybrid polytope in the $d$--dimensional lattice.  The first number, $k$,  is the $k$-dimensional simplicial element  $\sigma_k\in\mathcal S$, while the second starred number is its  co-dimension $k$  dual polytope in $\sigma^*_k\in{\mathcal S}^*$.   In the center column, we show the $k$--dimensional simplicial element in the black line, and its dual $(d-k)$--dimensional element in ${\mathcal S}^*$  is shown in red (grey). }\label{Fig:hybridvol}
\end{figure}

The circumcentric dual lattice is fundamental to RC \cite{FL:1984,M:1986,Caselle:1989,MM:2008b} which is not the case for dual lattices of other kinds, e.g. the barycentric dual lattice. Consider a  compact simplicial lattice ${\mathcal S}$ and its circumcentric dual  ${\mathcal S}^*$.  For simplicity, and not out of necessity, we  consider only  compact $d$--dimensional simplicial geometries in this paper.  This could be extended to geometries with boundaries by utilizing  previous results in RC \cite{Hartle:1981}.  Furthermore, we also assume here, for the same reasons,  that ${\mathcal S}$ is a well--centered Delaunay lattice \cite{Van:2009} and ${\mathcal S}^*$ is its Voronoi lattice.  The vertices, edges, triangles, $\ldots$ and $d$--simplexes in ${\mathcal S}$ will be denoted by $\sigma_0$, $\sigma_1$, $\sigma_2$, ..., $\sigma_d$; respectively.  The geometry of each simplex is determined uniquely by its $d \choose 2$ squared lengths; consequently the geometry of the entire simplicial lattice, $\left\{\ell_i^2\right\} \in {\mathcal S}$ is a function of all the edge lengths.  Similarly,  we label the vertices, edges, polyhedra and  polytopes of  ${\mathcal S}^*$ by $\sigma^*_d$, $\sigma^*_{d-1}$, $\sigma^*_{d-2}$, $\ldots$, $\sigma^*_1$, $\sigma^*_0$, respectively.  The geometry of each of these dual cells owes its existence, and can be expressed as,  functions of the squared edge lengths of $\mathcal S$.  The subscript on the simplicial elements is its dimension, e.g. $\sigma_k\in{\mathcal S}$ is a $k$-simplex;  where the subscript on the circumcentric dual polytope is its
co-dimension, e.g. $\sigma^*_k\in {\mathcal S}^*$ is the ($d-k$)-dimensional circumcentric dual element of $\sigma_k\in {\mathcal S}$.    Perhaps the most important property of the Voronoi and Delaunay lattices is that each element  $\sigma_k\in{\mathcal S}$ is perpendicular to its dual  element  $\sigma^*_k\in\mathcal{S}^*$.  The $d$--volume formed by connecting the vertices of  the $k$-simplex, $\sigma_k$, to the vertices of the ($d-k$)--dimensional polytope forms what we refer to as a $ d $--dimensional hybrid polytope, $V_{kk^*}$,  of volume,
\begin{equation}
\label{Eq:Vkk}
|V_{kk^*}| = \langle \sigma_{k}| \sigma_{k} \rangle = \int\, \sigma_k \wedge \sigma^*_k =  {d \choose k}^{-1}\, |\sigma_k| |\sigma^*_k| = \frac{k!(d-k)!}{d!}\, |\sigma_k| |\sigma^*_k|.
\end{equation} 
That this $d$--volume is a simple product of the volume of the dual--cell pairs reflects an important property of RC and is due to the orthogonality between the simplicial and circumcentric dual lattice.  In addition,  these $d$-hybrid volumes, when applied to all pairs of $k$-dual cells $\{\sigma^*_{k}, \sigma_{k}\}$, provide a proper tiling of the discrete geometry.  The curvature of the simplicial lattice has compact support on such  hybrid volumes, e.g. the $Rm$ is naturally defined on $V_{hh^*}$, the $Rc$ on $V_{\lambda \lambda^*}$, and the scalar curvature, $R$,  on $V_{vv^*}$, where $h$ is the $(d-2)$-dimensional hinge in ${\mathcal S}$, $\lambda$ is a dual edge in ${\mathcal S}^*$, and $v$ is a vertex in either ${\mathcal S}$ or ${\mathcal S}^*$  \cite{Hirani:2005,M:1997,MM:2008,AMM:2011}.  We provide an illustration of the five hybrid volumes in the case of a  4--dimensional lattice in Fig.~\ref{Fig:hybridvol}.

In constructing the curvature tensors on a simplicial geometry, one can take various approaches.  On a simplicial lattice, it is well known that the curvature is concentrated on the (d-2)--dimensional hinges $h$ \cite{Regge:1961}.  The form of this curvature is a conic singularity at the hinge. One approach used by researchers is to rigorously analyze this singularity structure, explicitly constructing a Dirac delta distribution for the curvature on the hinge $h$ by embedding the hybrid block and its simplicial support  into  a $(d+1)$-dimensional Euclidean space \cite{FL:1984,Sorkin:1975}.  This Dirac delta distribution is assumed to be uniformly distributed over the hinge, $ h $, and thus integral curvature measures were derived.  An alternative and equivalent  description has recently been used \cite{McDonald:2012,M:1986,M:1997,MM:2008,AMM:2011,MM:2008b}.  First, we take the support for this  curvature to be uniformly distributed over the entire hybrid block.  Second, we use the Voronoi area, $h^*$ and only this area to define the support for the sectional curvature.  Finally,  we assume the geometry within each hybrid block, $V_{hh^*}$, is an Einstein space  \cite{Eisenhart:1997}.  These approaches yield curvature components within each  hybrid cell that are proportional to each other, and can be transformed consistently to each other by taking suitable traces or identifying certain eigenvalues.  

In the remainder of this paper we assume that the components of curvature within the $V_{(d-2) (d-2)^*}$ hybrid block can be treated as an Einstein space, and we define the Gaussian curvature by the parallel transport of vectors restricted to the perimeters of  the dual Voronoi polygons, $h^*=\sigma^*_{d-2}$ where $h^*$   is the dual polygon to the ($d$-2)--dimensional simplicial  hinge, $h$. In this sense, the area of the dual Voronoi polygon,  $h^*$,  becomes the unique weighting area associated to hinge $h$ in order to define the Gaussian curvature, 
\begin{equation}
\label{Eq:GauCur}
\left( \begin{array}{c}Gaussian\ curvature \\  associated\ to\ hinge\ h\end{array}\right) = K_h = \frac{Angle\  Vector\  Rotates}{Area\  Circumnavigated} := \frac{\epsilon_h}{h^*},
\end{equation}
where $\epsilon_h$ is the deficit angle associated to hinge $h$, and $h^*$ is the area of the Voronoi plane dual to hinge $h$.  This  $V_{(d-2)(d-2)^*}$ hybrid cell is defined as the set of points in the lattice closer to the hinge $h$ than to any other hinge and provides a reasonable support for this curvature. In the next three subsections we briefly reexamine the construction of the Riemann ($Rm$), Ricci ($Rc$) and scalar curvature ($R$) on a simplicial lattice beginning as a guide with their respective definitions from the continuum. 

 For notational purposes, whenever we write a volume using a capital $V$, it will be understood that this represents a $d$--dimensional volume.  It can represent a $d$-dimensional hybrid or reduced--hybrid volume.  Its subscript will identify to which volume it corresponds.   For example, the hybrid cell $V_{(d-2)(d-2)^*}$ associated to hinge $h$ will be  written as $V_h=V_{(d-2)(d-2)^*}$. On the other hand, if we write, simply, $h$, in the correct context with no capitol, $V$, then it will mean the volume of $h=|h|$.  Our notation follows similarly for all other elements, e.g. $\ell=|\ell|$ for the length of an edge, $h^*= |h^*|$ for the area of the Voronoi polygon as in Eq.~\ref{Eq:GauCur}.

\subsection{The Simplicial {\em Rm}}
\label{subsec:Rm} 
In this section we examine the assumption alluded to in the previous section that in RC the geometry's  locally an Einstein space  within each hybrid cell, and  that the curvature is uniformly distributed within this block. To this end, one way we can construct the $Rm$ on $\mathcal S$ is  by examining the parallel transport of a vector around a closed loop.  In particular, a vector, $\pmb{v} := v^a\, {\pmb e}_b$ will ordinarily return rotated by an amount, $\pmb{e}_a\, \delta v^a$,  when parallel transported around the boundary of the area element, $A^{cd}\, \pmb{e}_c \wedge  \pmb{e}_d$,    
\begin{equation} 
\label{Eq:Riec}
\delta v^a = - 2\, R^a\phantom{}_{bcd}\, v^b\,  A^{cd}.
\end{equation}
This defines a constructive way to explore the components of the $Rm$ through the rotation operator.   
\begin{equation}
\label{Eq:ro}
\pmb{{\mathcal R}} =\frac{1}{4}\  \underbrace{{\pmb e}_a\wedge{\pmb e}_b}_{\begin{array}{c}  rotation\\ bivector\end{array}}\, 
\underbrace{R^{ab}\phantom{}_{cd}}_{Rm}\ \underbrace{{\pmb e}^c\wedge{\pmb e}^d}_{\begin{array}{c} loop\ of\\ circum\hbox{-}\\ navigation\end{array}}.
\end{equation}

The curvature of a $d$-dimensional piecewise flat simplicial manifold is concentrated on its co-dimensional 2 simplicial hinges, $h$ \cite{Regge:1961}.   Given any hinge, $h=\sigma_{d-2}\in\mathcal S$, there will be at least three $d$-simplicies sharing hinge $h$.  We assume here that there are $n\ge 3$ of these hinging on $h$, and we label these simplicies $\{ \sigma_d(i) \}$ for $i=1,2, \ldots, n$.  Let us consider a loop, $\mathcal C$,  within the interior of the collection of these $n$  $d$-simplicies that  encircles $h$ only once.  In addition, we consider a vector $\pmb{v}$ at a point $p_0$ along $\mathcal C$.  We place $p_0$ at the intersection of  $\mathcal C$ and  the ($d$-1)-simplex , $\sigma_{d-1}(1,2)$ that is shared by the two adjacent $d$-simplicies, $\sigma_d(1)$ and $\sigma_d(2)$.  The parallel transport of this vector within  $\sigma_d(1)$ is trivial since the geometry is flat.  The angle $\pmb v$ makes with the $d(d+1)/2$ edges of $\sigma_d$ remains constant as it is parallel transported along $\mathcal C$.  Furthermore, as we continue to parallel  transport  the vector along $\mathcal C$ into $\sigma_d(2)$ it is equally trivial. There is no curvature concentrated at their juncture ($\sigma_{d-1}(1,2)$), i.e. the two adjacent $d$-simplicies can be embedded into the same flat Euclidean geometry. We can continue to parallel transport $\pmb v$ along $\mathcal C$  through the remaining $(n-2)$ $d$-simplicies trivially; however, when we return to the starting point $p_0$ on $\mathcal C$ we will find the vector rotated by the deficit angle $\epsilon_h$, 
\begin{equation}
\epsilon_h := 2\pi - \sum_{i=1}^n  \theta_{i},
\end{equation}
where the summation is over all $\sigma_d(i)$ sharing hinge $h$, and  $\theta_{i}$ is the hyperdihedral angle between the two ($d-1$)-dimensional simplicial faces $\sigma_{d-1}(i-1,i)$ and $\sigma_{d-1}(i,i+1)$ of simplex  $\sigma_d(i)$ sharing hinge $h$.  {\em Remarkably, this rotation bivector lies in the plane, $h^*=\sigma^*_2$ perpendicular to hinge $h$, and this rotation is independent of the orientation and area of $\mathcal C$.}  This captures the integrated curvature at hinge $h$.  The curvature is a conic singularity.  The closed curve $\mathcal C$ can be shrunk to an arbitrarily small area, $A_{\mathcal C}$,  around $h$ without affecting the rotation in $\pmb v$ under parallel transport.  Therefore the sectional curvature can be made as large as one pleases by shrinking the area, 
\begin{equation}
\label{Eq:conic}
\left(\begin{array}{c} Conical\ Sectional\\ Curvature\ of\ hinge\\ h\in{\mathcal S}\end{array}   \right)  =\left(\frac{Angle\ vector\ rotates}{Area\ circumnavigated}  \right) = \lim_{A_{\mathcal C} \rightarrow 0} \  \frac{\epsilon_h}{A_{\mathcal C}}= \infty. 
\end{equation}

Is there a natural area that we can assign to the hinge, $h$?  What area does RC use?  To answer these questions, we identify three features of  $\mathcal S$ that yield a unique curvature measure.  First,  the simplicial geometry is piecewise flat.  Second,  each element  of simplicial geometry,  $\sigma_k \in \mathcal S$ is orthogonal to its circumcentric dual element, $\sigma_k^* \in {\mathcal S}^*$.   Thirdly, it has been shown that the Regge equations as well as the Regge action are uniquely defined as distributed quantities over such  circumcentric dual elements \cite{M:1986,M:1997,Lee:1984}. In his original paper on RC, Regge demonstrated that the curvature of a simplicial lattice is concentrated at each of the co-dimension 2 hinges, $h=\sigma_{d-2}\in{\mathcal S}$,  and this curvature  at each $h$ is a conic singularity as illustrated by Eq.~\ref{Eq:conic}.  However, in RC it appears that this curvature is defined to be distributed uniformly over the convex hull of the triangle hinge $h$ and the polygon, $h^*=\sigma^*_2\in {\mathcal S}^*$ dual to $h$.  Following the notation introduced in the end of last section, this area, $h^*\in{\mathcal S}^*$ is central to the hybrid cell $V_h:=V_{2 2^*}$, and this hybrid cell is the set of points in $\mathcal S$ closest to hinge $h$ than to any other hinge.  The area circumnavigated is,  $h^*:=|h^*|$ since RC defines the loop $\mathcal C$ to be the perimeter of the dual polygon $h^*$.  This area, $h^*$ and the deficit angle, $\epsilon_h$, are both  functions of the squared edges, $\{\ell^2_i\}$ of $\mathcal S$.   This (1) dual area $h^* \equiv |h^*|$, and (2) rotation, $\epsilon_h$, yield the hinge-based (or equivalently dual hinge-based) sectional curvature in RC,
\begin{equation}
\label{Eq:Kh}
K_h = K_{h^*}= \frac{\epsilon_h}{h^*}.
\end{equation} 

In order to use the continuum equations, Eqs.~\ref{Eq:Riec}-\ref{Eq:ro} we need to construct an appropriate set of basis vectors within the $V_{h}$ hybrid cell defined by the convex hull of $h$ and $h^*$.  Following the development in Sec.~\ref{sec:rf},
 we will find it convenient to define a time independent $d$-dimensional frame of  orthogonal bases 1-forms, $\{\pmb{e}^a\}$ for $a=1,2, \ldots, d$,  and  corresponding basis vectors, $\{\pmb{e}_a\}$ for $a=1,2, \ldots, d$, that  span a given hybrid polytope $V_{hh^*}$,
\begin{equation}
\pmb{e}_a = g_{ab}\, \pmb{e}^b.
\end{equation}
We refer to the following basis as the $\ell\lambda$--hybrid basis. We first identify an edge $\ell\in h$ and a dual edge $\lambda\in h^*$ in $V_{h}$.  If we connect the vertices of this $\lambda$ to the $d$ vertices of the hinge $h$, this will form a reduced hybrid $d$-simplex that we will refer to as $V_{h\lambda}$, i.e. the hybrid cell $V_h$ restricted to $\lambda$. The sum of all $V_{h\lambda}$'s for each $\lambda\in h^*$ equals  the volume of the non-simplicial hybrid cell,
\begin{equation}
V_h = \sum_{\lambda\in h^*} \, V_{h\lambda}.
\end{equation}
We define our orthogonal basis on each of these reduced hybrid cells, $V_{h\lambda}$.
In order to capture the orthogonality between $\mathcal S$ and ${\mathcal S}^*$ within the hybrid cell $V_{h}$, we define basis vectors of four kinds.
\begin{equation}
\label{Eq:elllambdabasis}
\pmb{e}_{a} = \left\{
\begin{array}{lll}
{\pmb e}_\ell  &=   \overrightarrow{\ell}, &\ell \in h,\\ 
{\pmb e}_i & =  \overrightarrow{m_{i}} = \overrightarrow{{\mathcal O}_i {\mathcal O}_{i+1}},  &  \forall i\in\{1, 2, \ldots,d-3\},\\
{\pmb e}_\lambda  &= \overrightarrow{\lambda},  &\lambda\in h^*,\\
{\pmb e}_{h \lambda}  &= \overrightarrow{m_{h \lambda}},
\end{array}
\right.
\end{equation} 
where, 
\begin{itemize} 
\item ${\mathcal O}_{1}$ is the circumcenter of edge $\ell \in h$,
\item ${\mathcal O}_{2}$ is the circumcenter of any one of the triangles, $\sigma_2$,  in $h$ such that $\sigma_2 \ni \ell$,
\item ${\mathcal O}_{3}$ is the circumcenter of any one of the tetrahedra in $h$ such that $\sigma_3 \ni \sigma_2$,
\item[\, ] \centerline{$\vdots$} 
\item ${\mathcal O}_{d-3}$ is the circumcenter of one of two the (d-3)-dimensional simplexes in $h$ such that  $\sigma_{d-3} \ni \sigma_{d-4}$, 
\item ${\mathcal O}_{d-2}$ is the circumcenter of the $(d-2)$-dimensional simplicial hinge, $h = \sigma_{d-2}$,
\end{itemize}
and, $m_{\ell \lambda}$ is the segment from ${\mathcal O}_1$ to the point on lambda, ${\mathcal O}_\lambda^*$,  so that this line segment is  perpendicular  to  $\lambda$, i.e. it intersects $\lambda$ at the circumcenter of the ($d$-1)-dimensional simplex, $\lambda^*$.  The set of (d-2) mutually orthogonal vectors,\footnote{In 2-dimensions there  the hinge is a vertex and the space is a panned by only two orthogonal vectors, $\{\vec \lambda,\vec m_{h\lambda}\}$, i.e. there are no $\vec m_i$'s.  In three dimensions there are also no $m_i$'s and the space is a panned by the three orthogonal vectors, $\{\vec \ell, \vec \lambda,\vec m_{h\lambda}\}$.  It is only in four and higher dimensions that we start generating the $m_i$'s.} 
\begin{equation}
\left\{ \overrightarrow{\ell}, \overrightarrow{m_{i}}\right\},
\end{equation}
 span the hinge, $h=\sigma_{d-2}$, while the remaining two vectors, 
 \begin{equation}
 \label{Eq:hstar}
 \left\{ \overrightarrow{\lambda},\overrightarrow{m_{h \lambda}} \right\}
 \end{equation}
  span the dual polyhedron, $h^*$.   The hinge can be expressed as a ($d$-2)-vector, 
\begin{equation}
\label{Eq:hinge}
{\pmb h} = \frac{1}{(d-2)!}\left( \pmb{e}_{\ell}\wedge\pmb{e}_{1}\wedge  \pmb{e}_2 \wedge  \ldots \wedge \pmb{e}_{d-3} \right),
\end{equation}
where the ($d$-2)-volume of the hinge is $h = |{\pmb h}|$.  In addition, the dual hinge can be equivalently expressed either as the dual of Eq.~\ref{Eq:hinge}
or by  using the last two basis vectors in Eq.~\ref{Eq:hstar},
\begin{equation}
\label{Eq:hstarhat}
\pmb{h}^* = h^*\, \frac{  \pmb{e}_{\lambda}\wedge\ \pmb{e}_{h\lambda} }{ |\pmb{e}_{\lambda}\wedge\ \pmb{e}_{h\lambda} |} = h^* \hat{\pmb h}^*.
\end{equation}
The contravariant components of this flat-space metric, $\pmb g$, in this vector basis is,
\begin{equation}
g^{ab} = diag\left(\underbrace{\ell^{-2},  m_0^{-2}, m_1^{-2}, \ldots, m_{d-3}^{-2}}_{hinge,\ h}, \underbrace{\lambda^{-2}, m_{h\lambda}^{-2}}_{h^*} \right).
\end{equation}
A dual set of orthogonal basis 1-forms can be generated, 
\begin{equation}
{\pmb \omega}^a = g^{ab}\, {\pmb e}_b,
\end{equation}
and the covariant components of the metric are just, 
\begin{equation}
\label{Eq:metricdown}
g_{ab} = diag\left(\underbrace{\ell^2, m_0^2, m_1^2, \ldots, m_{d-3}^2}_{hinge, \ h}, \underbrace{\lambda^2, m_{h\lambda}^2}_{h^*}\right).
\end{equation}

We are now in a position to construct the $Rm_h$ associated to hinge $h$ by examining the rotation operator, $\pmb{\mathcal R}$, defined by  Eq.~\ref{Eq:ro}. The sectional curvature, $K_h$, is naturally defined on the convex hull of $h$ and $h^*$, where the convex hull is just the hybrid volume,  $V_{h}=V_{22^*}$ (as  illustrated in Fig.~\ref{Fig:hybridvol} for the special case of 4-dimensions).  Furthermore, this is the only nonzero sectional curvature within the hybrid polytope $V_{h}$.    It is  the orthogonality between any element of  the simplicial lattice, $\sigma_k$,  and its corresponding circumcentric dual, $\sigma^*_k$,  that yields an Einstein space geometry within each hybrid cell, and that each of the $d(d-1)$ non--zero components of the $Rm_h$ is proportional to the sectional curvature $K_h$\cite{M:1997}.  No matter what the orientation of the loop, $\pmb \Sigma$,  we take encircling the hinge $h$, the unit rotation bivector will always lie in the plane of $h^*$ and have magnitude $\epsilon_h$,  
\begin{equation} 
\label{Eq:firsttwo}
\left(\ .\ \right)\pmb{{\mathcal R}}\left(\ \pmb{\Sigma}\ \right)  =  \pmb{\hat{h}}^*\  \epsilon_h.
\end{equation} 
Therefore, the only non-vanishing components of the simplicial $Rm$ associated to hinge $h$ are when the first two indices lie in the dual polygon, $h^*$. While, it is true that in the basis chosen above the only non vanishing components of the $Rm$ are in the Voronoi plane, $h^*$,
\begin{equation}
R^{\lambda m_\lambda}{}_{\lambda m_\lambda} =  \frac{\epsilon_h}{h^*}.
\end{equation}
Nevertheless, this would not describe adequately the curvature within the hybrid cell. In particular, any basis chosen in general position would ordinarily (except for a set of measure zero) not have a single vector that lies in the hyperplane spanned by  $h$.  This would yield a fully--populated rotation matrix.   In particular, this would be true for our specially constructed  basis if we rotated each of the basis vectors $\left\{ \vec \ell, \vec m_1, \vec m_2, \ldots \vec m_d-3\right\}$ even slightly out of the plane of $h$.  {\em In this sense, the space within each hybrid cell, $V_h$ , is a Einstein space} \cite{Eisenhart:1997}. Additionally, if we restrict the parallel transport of vectors around the perimeter of the Voronoi polygon, $h^*$, then we have usual Einstein space expressions for those components of the Riemann, $Rm_h$, the Ricci, $Rc_h$ and scalar $R_h$ curvatures associated exclusively to the Einstein space of each hybrid cell $V_{h}$ ,
\begin{eqnarray}
\label{Eq:Rmh}
Rm^{(hyb)}_h & = &   K_h = \frac{\epsilon_h}{h^*},\\
\label{Eq:Rch}
Rc^{(hyb)}_h & = & (d-1)\, K_h  =  (d-1)\,  \frac{\epsilon_h}{h^*}, \ \ \hbox{and} \\
\label{Eq:Rh}
R^{(hyb)}_h & = & d(d-1)\, K_h  =  d(d-1)\,  \frac{\epsilon_h}{h^*}.
\end{eqnarray}
Each of these Einstein space curvatures (Eqs.~\ref{Eq:Rmh}-\ref{Eq:Rh}) associated with their hybrid cell are uniformly distributed within this cell. This interpretation is consistent with other derivations \cite{FL:1984,Sorkin:1975} with proper care in their  normalization.  The component of the Riemann associated to the hinge-based hybrid cell is the entire Riemann, 
\begin{equation}
Rm^{(hyb)}_h = Rm_h,
\end{equation}
as it realizes the entire sectional curvature attributable to the hinge $h$.  Therefore,  the Riemann tensor in Cartan language is,     
\begin{equation}
\label{Eq:Rm}
\pmb{Rm}_h =  \pmb{\hat h}^*\,  Rm_h\  \pmb{\hat h}^*=  \pmb{\hat h}^*\  K_h\ \pmb{\hat h}^* = \pmb{\hat h}^*\ \frac{\epsilon_h}{h^*}\ \pmb{\hat h}^*, 
\end{equation}
where,  $\pmb{\hat h^*} =  \pmb{h}^*/h^*$ is the unit 2-form generated  on $h^*$ in hybrid volume, $V_{h}$.  In the $\ell\lambda$--hybrid basis this unit bivector  is given by  Eq.~\ref{Eq:hstarhat}.

As an illustrative example, we can use these expressions for the $Rm_h$ and the hybrid proper volume,  $V_{h}$,   to derive the RC expression for the Hilbert action \cite{Regge:1961,M:1997}.  In RC the geometric expressions are often expressed as integrated quantities over the relevant d-dimensional hybrid or reduced-hybrid volume elements, Fig.~\ref{Fig:hybridvol}.  The Regge-Hilbert action, $I_{R\hbox{-}H}$,  in general relativity, is one such integrated quantity.    
\begin{eqnarray} 
I = \frac{1}{16\pi}\,\int R\,dV_{proper} & \approx &  \frac{1}{16\pi}\, \sum_{h}\  R^{(hyb)}_h\, V_{h}\\ 
& = & \frac{1}{16\pi}\, \sum_{h}\   \left(d(d-1)\, \frac{\epsilon_h}{|h^*|}\right)\, \left(  \frac{2}{d(d-1)}\, |h|\, |h^*|\right)\\ 
& = &    \frac{1}{8\pi}\,  \sum_{h}\  \epsilon_h\,h.
\end{eqnarray}

\subsection{The Simplicial Ricci Tensor ({\em Rc}) }
\label{subsec:Rc} 
In the last section we constructed the Riemann tensor and in this section we  examine the simplicial $Rc$ on $\mathcal S$.  This will involve the simplicial version of taking the trace of the simplicial $Rm$ tensor in Eq.~\ref{Eq:Rmh}. We examine the $Rc$ tensor here as it appears on the  right hand side of the RF equation (Eq.~\ref{Eq:RF}).   This tensor was first  introduced in \cite{FL:1984},  and later in context within the hybrid cells by \cite{AMM:2011}.  We briefly reconstruct this tensor by utilizing the  Einstein space geometry interior to $V_h$,  and by using the $\ell\lambda$-hybrid basis construction (Eq.~\ref{Eq:elllambdabasis})  discussed in the previous section.  This construction provides us with a  geometric description of the trace operation on a lattice geometry.  

\begin{figure}
\centering
\includegraphics[height=3.0in]{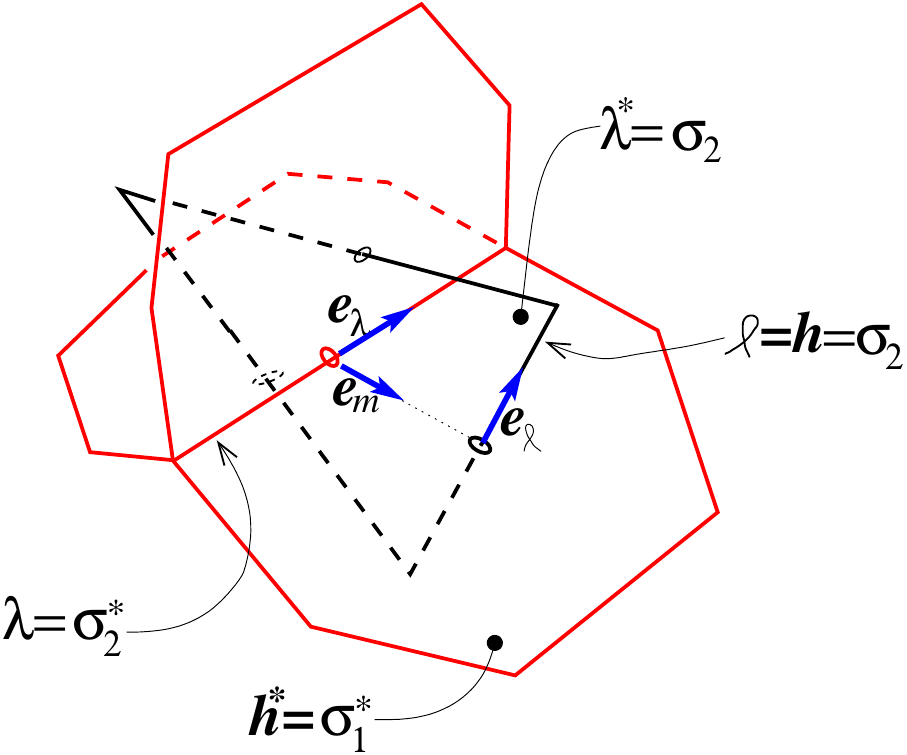}
\caption{We illustrate here the lattice geometry used to define the $Rc$ tensor for a 3-dimensional simplicial lattice.  The $Rc$ tensor is naturally defined on an edge, $\lambda=\sigma^*_1\in \mathcal S^*$.  Along this edge we define one of the unit triad vectors, $\hat {\pmb e}_\lambda$ of the $\ell\lambda$-basis, where ${\pmb e}_\lambda = \lambda\, \hat {\pmb e}_\lambda$ as usual.  In this illustration, we show the three dual polygons in red, $\ell^*_{|\lambda}\in \mathcal S^*$ sharing edge $\lambda$.  Dual to each of these polygons is an edge of the triangle $h=\lambda^*\in \mathcal S$.  Along one of these edges, $\ell \subset h$,  we define the second of the triad vectors, ${\pmb e}_\ell = \ell\, \hat{\pmb e}_\ell$.  Finally the normalized vector from the center of edge $\ell$ perpendicular to edge $\lambda$, defines our final triad vector, $\pmb{e}_{\ell\lambda} = m_{\ell\lambda}\, \hat{\pmb e}_m$.}\label{Fig:Rclambda}
\end{figure}

In the continuum, the simplicial $Rc$ is a 1-form valued 1-form, 
\begin{equation}
\label{Eq:RcRm}
\pmb{Rc} = {\pmb e}^a\, R_{ab}\, {\pmb e}^b = {\pmb e}^a\, R^c\phantom{}_{acb}\, {\pmb e}^b.
\end{equation}
This tensor is defined along some basis vector, ${\pmb e_b}$ in the tangent space.  This basis vector defines a common set of area bivectors for circumnavigation, i.e. ${\pmb \Sigma}=\sigma^{cb} {\pmb e}_c \wedge {\pmb e}_b$, one for each, $c\ne b$ .  In the previous section we defined the $Rm_h$ in RC on the dual Voronoi areas,  $h^*=\sigma^*_{d-2}\in{\mathcal S}^*$ of the hybrid cell, $V_{h}$.   Each edge, $\lambda=\sigma^*_{d-1}\in{\mathcal S}^*$, is a common bounding edge of $d(d-1)/2$ of these dual polygons, $h^*$.  It is therefore natural to define a simplicial Ricci tensor, $Rc_\lambda$, along these dual edges, $\lambda$.  In this sense the trace of the $Rm$ is a suitable sum over each of these Voronoi ares, $h^*$ sharing edge $\lambda$.  We illustrate this in Fig.~\ref{Fig:Rclambda} for the case of three dimensions.  

\begin{figure}
\centering
\includegraphics[height=2.0in]{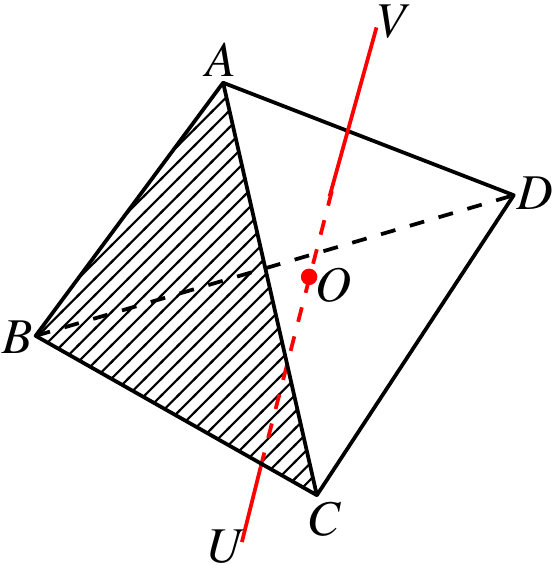}
\caption{This figure illustrates the various geometric elements that were used to define the $Rc_\lambda$ in 4-dimensions.  The line segment $\lambda = \overline{UV}$ is the edge in the circumcentric dual lattice, $\mathcal S^*$, where we  define the Ricci tensor. The tetrahedron $\lambda^*=\{ABCD\}\in {\mathcal S}$ is perpendicular to $\lambda$.  If we connect the endpoints of $\lambda$ to the four vertices of $\lambda^*$ the resulting 4-dimensional cell is the $\lambda$-hybrid simplex of volume, $V_\lambda = (1/d) \lambda \lambda^*$.  In 4-dimensions the hinges are triangles.  We highlight one of the 4 triangles, $h=\{A,B,C\}\subset \partial \lambda^*$, bounding $\lambda^*$ by shading the triangle.   The point $O$ lies at the circumcenter of $\lambda^*$.  It is the point equidistant from vertices, $A$, $B$, $C$ and $D$.  Point $O$ also lies on the dual edge $\lambda$.  We refer to the tetrahedron formed by connecting $O$ with triangle hinge, $h=\{A,B,C\}$,  as the restricted  3-volume, $\lambda^*_h$.  There are four such reduced 3-volumes subdividing,  $\lambda^*$.  If we connect this shaded hinge, $_h$, to edge $\lambda$, the resulting simplex,  will be the $h\lambda$-reduced hybrid block of volume,  $V_{\lambda h} = (1/d) \lambda h_\lambda$.   }\label{Fig:Vlamstarh}
\end{figure}

We construct the Ricci tensor, $Rc_\lambda$ ,  associated to a dual edge  $\lambda$ by  properly tracing over each of the Ricci tensors, $Rc^{(hyb)}_h$, defined in association with each Voronoi polygon,  $h^*$, containing dual edge $\lambda$.  The trace, or sum, will be over  all  $h^*$ restricted to $\lambda$, or $h^*{}_{|\lambda}$. The trace  is accomplished using a weighted average as described for RC in  \cite{M:1997,MM:2008,AMM:2011} and independently in discrete exterior calculus by \cite{Hirani:2005}.  Here one considers the component of the Ricci tensor associated to the hybrid cell, $Rc_\lambda$, to be uniformly distributed over the  hybrid $d$-volume of the simplicial lattice associated with edge $V_\lambda$ and integrates this over this volume.  This integrated tensor can be identified  with the trace of the familiar hinge-based expression for simplicial Riemann tensor, $Rm_h$  in Eq.~\ref{Eq:Rch} with the inclusion of a proper weighting factor.  In particular, 
\begin{equation}\label{Eq:lambdaRic}
Rc_\lambda \, V_\lambda = \underbrace{\sum_{h_{|\lambda^*}} Rc^{(hyb)}_{h} \, V_{\lambda h}}_{\begin{array}{c} simplicial\,form\\of\,the\,trace\end{array}},
\end{equation}
where, 
\begin{eqnarray}
V_\lambda &=&  {d  \choose 1}^{-1}\, \lambda\  \lambda^* = \sum_{h_{|\lambda^*}} {d  \choose 2}^{-1}\, h\, h^*_\lambda,\\
V_{\lambda h} &=&  {d  \choose 2}^{-1}\,  h\, h^*_\lambda, \ \ \ \hbox{and}\\
Rc^{(hyb)}_h &=&  (d-1) Rm_h = (d-1)\, \frac{\epsilon_h}{h^*} .
\end{eqnarray}
Hence, 
\begin{equation} 
\label{Eq:RCLambda}
Rc_\lambda = \frac{\sum_{h_{|\lambda^*}} Rc^{(hyb)}_h\, h h^*_\lambda}{\sum_{h_{|\lambda^*}}  h h^*_\lambda} :=  \langle Rc^{(hyb)}_h \rangle_\lambda .
\end{equation}
In constructing $Rc_\lambda$ we make use of reduced hybrid blocks, $V_{h^* \lambda}$, as well as reduced hinge volumes, $h^*_\lambda$ which are illustrated in the case of $d=4$ in Fig.~\ref{Fig:Vlamstarh}.  In general, the reduced hybrid $d$-volume, 
$V_{h^* \lambda}$, is the fraction of the set of points interior to the hybrid $d$-volume, $V_{h^*}$ that are closest to the dual edge, $\lambda$.  This provides a decomposition of the hybrid cell into reduced hybrid cells, 
\begin{equation}
V_h = \sum_{\lambda\in h^*} V_{h\lambda}, 
\end{equation}
as illustrated in Fig.~\ref{Fig:Rclambda}.
Each reduced hybrid cell,  $V_{h\lambda}$, contain the circumcenters, ${\mathcal O}_i$,  that defined the $\ell\lambda$--hybrid basis.  Similarly, the reduced hinge is a ($d$-2)-volume, $\lambda^*_h$,  and is defined to be the set of points in the interior of the ($d$-1)-dimensional simplex, $\sigma_\lambda^*$,  that are closest to its bounding  hinge, $h$, than to any other hinge in its boundary.

The last two indices, $c$ and $b$, of the $Rm$ in the  far right-hand side of Eq.~\ref{Eq:RcRm} identify the components of a  circumnavigated area, and in our case this will be the Voronoi polygons, $h^*_{|\lambda}$, sharing edge $\lambda$.  In the $h\lambda$-basis the two indices identifying this area element are  $b=\lambda$ and $c=m_{h\lambda}$, or vise versa.  As we argued earlier, the $Rc$ tensor is associated with dual edge, $\lambda$,  consequently $b=\lambda$ and $c=m_{h\lambda}$. The Ricci tensor is diagonal in the space of bivector-valued 2-forms.  Therefore, since index $b=\lambda$ then   index $a$ must be the moment arm $m_{\ell\lambda}$.  There is a unique $m_{h\lambda}$ for each of the dual $d \choose 2$ Voronoi polygons, $h^*_{|\lambda}$, sharing edge $\lambda$.  We define $Rc_\lambda$ by summing over these $h\lambda$-moment arms. Therefore, we would sum over all all polygons, $h^*$, sharing edge, $\lambda$, i.e. 
\begin{equation}
\label{Eq:Rcfirst}
{\pmb Rc}\left( \ldots, \pmb \lambda \right) = {\pmb e}_a\, R^{ca}\phantom{}_{cb}\, {\pmb e}^b\left( {\pmb e}_\lambda \right)= 
{\pmb e}_a\, R^{ca}\phantom{}_{c\lambda} = \sum_{h^*_{|\lambda}}\,   {\pmb e}_a\, R^{m_\lambda a}\phantom{}_ { m_\lambda \lambda}. 
\end{equation}
Finally, since we showed that the simplicial $Rm$ has compact support and uniformly distributed over the dual hybrid polytopes, $V_{(d-2)(d-2)^*}$,  and since the geometry within this hybrid block is locally an Einstein space, then he only non-zero vector component, $a$  left in Eq.~\ref{Eq:Rcfirst} is along $\lambda$.  Therefore, the simplicial $Rc_\lambda$ is diagonal, 
\begin{equation}
{\pmb Rc}\left( \ldots, \pmb \lambda \right) = \sum_{h^*_{|\lambda}}\,   {\pmb e}_\lambda\, R^{m_\lambda \lambda}\phantom{}_ { m_\lambda \lambda} = \sum_{h^*_{|\lambda}}\,   {\pmb e}_\lambda\, R^{h^*}\phantom{}_ { h^*} = {\pmb e}_\lambda\,  \sum_{h^*_{|\lambda}}\,   R^{h^*}\phantom{}_ { h^*}. 
\end{equation}
  
\begin{figure}
\centering
\includegraphics[height=1.75in]{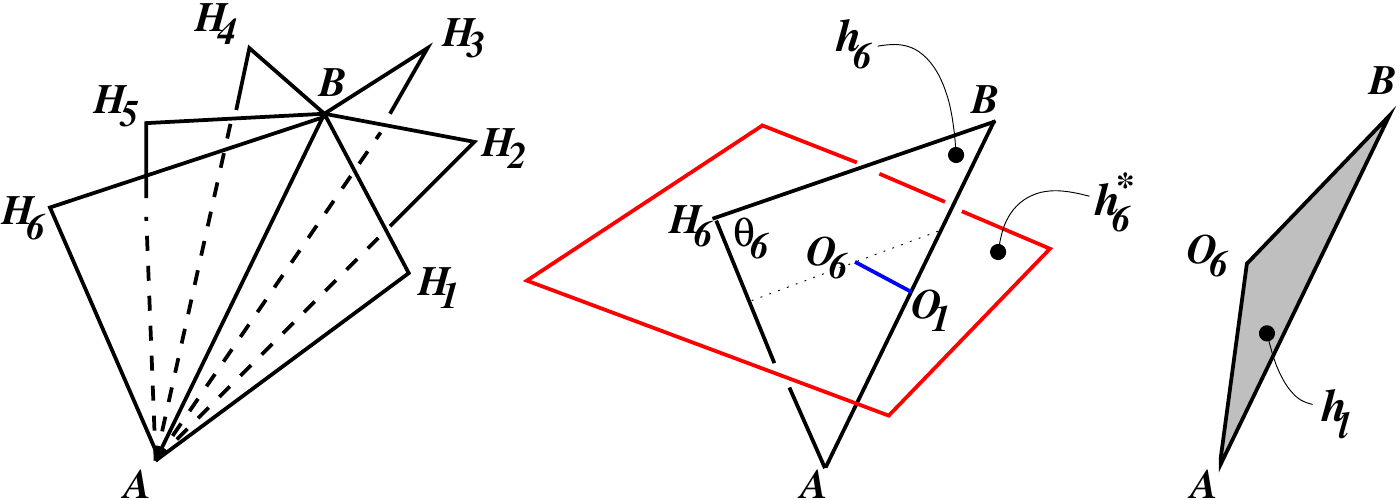}
\caption{In the left-most figure, we focus our attention on an edge $\ell=\overline{AB}$ in a 4-dimensional simplicial lattice, $ {\mathcal S}$, and  illustrate the hinges $h$ sharing this edge.  In this particular case there are six such hinges, $h_i = \overline{ABH_i}$, for $ i=1,2,\ldots,6$. The orientation on edge $\vec{\ell}=\vec{AB}$ induces an orientation on each of the six triangles sharing the edge.  In particular, the triangle is a bi-vector in ${\mathcal S}$, where $h_i= |h_i| \, \hat{h} = |ABH_i|\ \left( \frac{\vec{AB} \wedge \vec{AH_i}}{2 |ABN_i|} \right)$.  In the central figure we focus our attention on one of the six hinges, $h_6$,  sharing edge $\ell$, and its Voronoi dual polygon $h^*_6\in {\mathcal S}^*$.   This dual polygon, $h^*_6$  is perpendicular to the hinge, $h_6$. The dotted line shown in the figure is simply the intersection of these orthogonal 2--dimensional surfaces.  The orientation of the edge, $\ell$ also induces an orientation on the bi-vector, $h_6^* \in {\mathcal S}^*$.  In particular, $
{\pmb h}_i^* =  |h_i^*| \, \hat{\pmb h}_i^* =  h_i^* \, \left( \hat{\pmb h}_i\right)^*=  
h^*_i \ \left( \frac{\vec{AB} \wedge \vec{AH}_i}{2 |ABN_i|}\right)^*.$
The circumcenter of hinge $h_6$ is labeled  ${\mathcal O}_6$ .  This point is   equidistant from vertices $A$, $B$ and $H_6$.  The ``$h\lambda$--moment arm''  is the edge $m_{h\lambda}=\overline{O_1 O_6}=(1/2)\vec{\ell} \cot{(\theta)}$ reaching from the center of edge $\ell$ to the circumcenter of triangle $h_6$.  It  is instrumental in the Cartan construction of the Regge--Einstein tensor \cite{M:1986}. The right--most figure shows the fraction of hinge $h_6$ closer to edge $\ell$ than to its other two edges.  This edge--based triangle in ${\mathcal S}$ has area $h_\ell = \overline{ABO_6} = 1/2\, \ell^2 \cot{(\theta_6)}$, which is useful in defining the $Rc$ in RC.}
\label{Fig:RicEdge}
\end{figure}

 To complete the analysis of the right hand side of the RRF equation, we  need also to define the Ricci tensor, $Rc_\ell$, on an edge $\ell\in {\mathcal S}$.  In other words we need to express this $\ell$-based Ricci tensor associating it with an  appropriate weighted sum of the $\lambda$-based Ricci tensor  equation (Eq.~\ref{Eq:RCLambda}) onto the simplicial edge $\ell$.  In essence, we trace over the components of $Rc_\lambda$, one component for each $\lambda \in \partial\ell^*$,  with respect to the simplicial edge $\ell$  \cite{McDonald:2012}.  We will have to similarly trace over the left-hand side of the RRF equation in Sec.~\ref{sec:theorem} in order to prove the Corollary 1.  In particular, consider an edge $\ell \in {\mathcal S}$ as illustrated in Fig.~\ref{Fig:RicEdge}.  Similarly, one integrates the $Rc_\ell$ 
over the hybrid $d$-volume of the simplicial lattice associated with edge $V_\ell$ and equates this with the trace of the familiar hinge-based  simplicial Ricci tensor ($Rc_h$) \cite{M:1997},
\begin{equation}\label{Eq:IRic}
Rc_\ell \, V_\ell = \underbrace{\sum_{\lambda_{|\ell^*}} Rc_\lambda\, V_{\ell \lambda}}_{\begin{array}{c} simplicial\,form\\of\,the\,trace\end{array}},
\end{equation}
where the sum is over the dual edges, $\lambda$ bounding Voronoi polygon, $\ell^*$,  and 
\begin{eqnarray}
V_\ell &=&  {d  \choose 1}^{-1}\, \ell\  \ell^* = \sum_{h_{|\ell}} {d  \choose 2}^{-1}\, h_\ell\  h^*,\ \hbox{and}\\
V_{\ell \lambda} &=&  {V_\ell}_{|\lambda} = \sum_{h_{|{\ell}}} {d  \choose 2}^{-1}\, {h_\ell}_{|\lambda}\  h^*.
\end{eqnarray}
In order for us to obtain  the right-hand side of Eq.~\ref{Eq:srfell} we make use of  the following definition for the edge-based weighted average:
\begin{equation}
\langle f_h \rangle_\ell :=  \frac{1}{V_\ell} \sum_{h_{|\ell}} f_h\,V_{\ell h} = \frac{\sum_{h_{|\ell}}\, f_h h_\ell}{\sum_{h_{|\ell}}\,  h_\ell}.
\end{equation}
This yields an explicit expression for the Ricci tensor, $Rc_\ell$,  associated to an edge, $\ell$,  in the simplicial geometry,
\begin{eqnarray} 
Rc_\ell & =  & \frac{1}{V_\ell} \,\sum_{\lambda_{|\ell^*}} Rc_\lambda\, V_{\ell \lambda} \\
            & =  & \frac{1}{V_\ell} \,\sum_{\lambda_{|\ell^*}}  \left(  \frac{1}{V_\lambda} \,\sum_{h_{|\lambda^*}} Rc^{(hyb)}_{h}\, V_{\lambda h}  \right) V_{\ell\lambda} \\
            & = & \sum_{\lambda_{|\ell^*}} \sum_{h_{|\lambda^*}} \, Rc^{(hyb)}_{h}\, \left(\frac{V_{\lambda h}}{V_\lambda}\right) \, \left(\frac{V_{\ell\lambda}}{V_\ell}\right) ,
\end{eqnarray}
or equivalently, 
\begin{equation} \label{Eq:Rcell}
Rc_\ell = (d-1)\,  \langle \langle Rm_h \rangle_\lambda \rangle_\ell.
\end{equation}

\subsection{The Simplicial {\em R}}
\label{subsec:Rs} 
\begin{figure}
\centering
\includegraphics[height=4.0in]{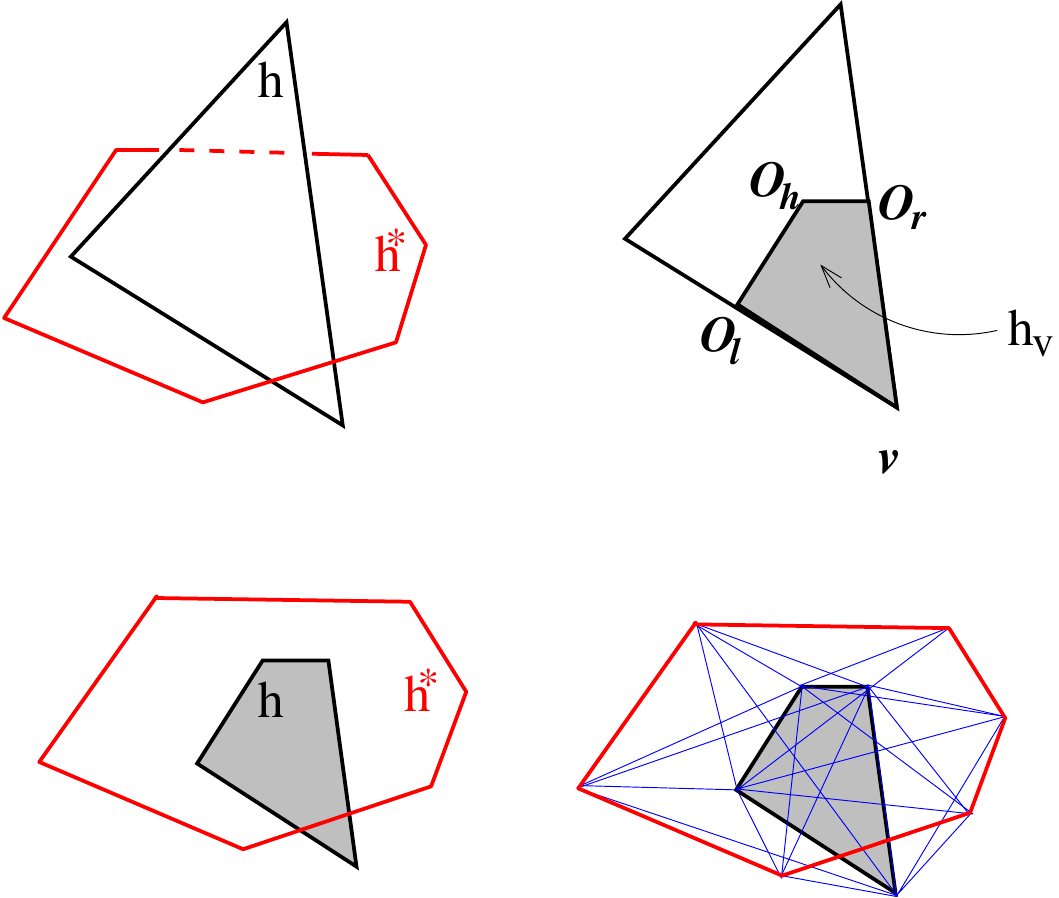}
\caption{  Here we show the reduced hybrid $h_v$ volume restricted to the vertex $v=\sigma_0 \in {\mathcal S}$. This reduced hybrid d--dimensional polytope, $h_v$, is formed by connecting the kite--shaped set of points in the triangle $h \in {\mathcal S}$ closest to the vertex $v$ to the dual area, $h^*\in {\mathcal S}^*$.
}\label{Fig:hv}
\end{figure}

The scalar curvature on the lattice is naturally defined at a vertex, $v\in{\mathcal S}$ \cite{MM:2008}. The construction of the scalar curvature, $R_v$,  on a simplicial lattice is similar to the derivation of $Rc_\ell$ in the last section, but less involved.  The scalar curvature requires a proper weighted average over the hinge--based Riemann curvatures $Rm_h$ of  Eq.~\ref{Eq:Rh}.  As we observed in Sec.~\ref{subsec:Rm}, the $Rm_h$ in RC is defined over  the hinge--based hybrid volumes, the Ricci tensor we derived in  Sec.~\ref{subsec:Rc} was naturally associated with an edge-based hybrid volume, $V_\ell$,  and the scalar curvature $R_v$ is naturally defined on the vertices, $v\in \mathcal S$, over a vertex-based reduced hybrid cell $V_v$.  This trace requires a sum over a suitable reduced hybrid volume \cite{MM:2008}, and in particular the volume, $V_{hv}$.  This reduced reduced hybrid volume is the set of points of the hinge-based hybrid cell, $V_h$ that are closest to one of its vertices, $v$ than to any other.  We illustrate this reduced hybrid cell for the case of $d=4$ in Fig.~\ref{Fig:hv}.
In particular we find, 
\begin{equation}\label{Eq:vRic}
R_v \, V_v = \underbrace{\sum_{h_{|v}} R^{(hyb)}_h\, V_{hv}}_{\begin{array}{c} simplicial\,form\\of\,the\,trace\end{array}},
\end{equation}
where, 
\begin{eqnarray}
V_v &=&  {d  \choose 0}^{-1}\,  v^* = \sum_{h_{|v}} {d  \choose 2}^{-1}\, h_v\  h^*,\\
V_{hv} &=& \sum_{h_{|v}}  {d  \choose 2}^{-1}\, h_v\  h^*,\ \hbox{and}\\
R^{(hyb)}_h &=&  d(d-1)\, \frac{\epsilon_h}{h^*}.
\end{eqnarray}
This yields an explicit expression for the Riemann scalar, $R_v$,  associated to a vertex, $v\in {\mathcal S}$,
\begin{eqnarray} 
R_v & = & \frac{1}{V_v} \sum_{h_{|v}} d(d-1) \frac{\epsilon_h}{h^*}\, V_{hv} \\
& = &  \frac{\sum_{h{}_{|v}} \epsilon_h\, h_{v}}{\sum_{h{}_{|v}}  h^*\, h_{v}} \\
& = &  d(d-1) \frac{\langle \epsilon_h \rangle_v}{\langle h^* \rangle_v},
\end{eqnarray}
or equivalently, 
\begin{equation}
R_v = d(d-1)\, \langle Rm_h \rangle_v. 
\end{equation}

\subsection{Overview: the {\em Rm}, {\em Rc} and {\em R}}
\label{subsec:RmRcRs} 

In the last three subsections, we constructed the Riemann, Ricci and scalar curvatures on a $d$--dimensional simplicial lattice, $\mathcal S$.  The coupling of elements of the Delaunay lattice to dual elements of  the Voronoi lattice appear fundamental in their construction.  The orthogonality between elements in $\sigma_k \in {\mathcal S}$ and their circumcentric duals, $\sigma_k^*\in{\mathcal S}^*$ gave rise to relatively simple diagonalized expressions for the curvature expressions.  This closed-form factorization of the homology and cohomology structures admits a unique interpretation of the curvatures in terms of hybrid blocks and reduced hybrid blocks.  These hybrid blocks are  the product of Delaunay and dual  Voronoi elements, 
\begin{equation}
V_k =  {d \choose k}^{-1}\, \sigma_k \wedge \sigma^*_k.
\end{equation}
We framed the curvatures by (1) uniformly distributing the curvature within their respective hybrid cells, (2) restricting the domain of parallel transport to the scaffolding provided by the dual circumcentric lattice, and (3) identifying the interior of each hybrid block with an Einstein space geometry, where the mixed Ricci tensor is simply
\begin{equation}
Rc^{(hyb)}_h = \frac{1}{d}\, R^{(hyb)}_h= (d-1)\, Rm_h.
\end{equation}
We showed in the last three sections that the part of simplicial Riemann, Ricci and scalar curvatures associated with the hybrid cells can be defined, most naturally, on the hinge--based, edge--based and vertex--based hybrid cells, respectively.  In particular, the expressions for the curvatures that we will use, in part,  to derive the right--hand side of the RRF equations in the next section, are as follows: 
\begin{eqnarray}
Rm_h & = & \langle Rm_h \rangle_h   = \frac{\epsilon_h}{h^*}\\
Rc_\ell & = &  (d-1)\, \langle\langle Rm_h \rangle_\lambda \rangle_\ell  \\
R_v & = & d(d-1)\,  \langle Rm_h \rangle_v ;
\end{eqnarray}
 and,
 \begin{equation} 
 Rc_\lambda =  \langle Rc_h \rangle_\lambda   =  (d-1)\, \frac{\sum_{h_{|\lambda^*}} h\, \epsilon_h\, \left(h^*_\lambda/h^* \right)}{\sum_{h_{|\lambda^*}} h\, h^*_\lambda}.
 \end{equation}
Additionally, we  can project or assign any one of these curvatures to the vertices, edges and hinges of $\mathcal S$ as shown in Table~\ref{tab:RmRcR}.
 \begin{table}[h] \begin{center}
\begin{tabular}{| c || c | c | c |}
\hline
{ } & hinge, $h$  & edge, $\ell$ & vertex, $v$ \\ \hline \hline
$Rm$ & $\boxed{\langle Rm_h \rangle_h = Rm_h}$&  $\langle Rm_h \rangle_\ell$ & $\langle Rm_h \rangle_v$\\ \hline
$Rc$ & $(d-1)\, \langle Rm_h \rangle = (d-1) Rm_h$ &  $\boxed{(d-1)\, \langle \langle Rm_h \rangle_\lambda \rangle_\ell}$ & $(d-1)\, \langle \langle Rm_h \rangle_\lambda \rangle_v $\\ \hline
$R$ & $d(d-1)\, \langle Rm_h \rangle = d(d-1) Rm_h$ &  $d(d-1)\, \langle Rm_h \rangle_\ell$ & $\boxed{d(d-1)\, \langle Rm_h \rangle_v}$\\ \hline
\end{tabular}
\caption{The expressions for the simplicial curvatures. The three boxed equations on the diagonal are particularly useful in the RRF equation construction.  The expression for the Riemann scalar will be useful if we consider normalized Ricci flow.}\label{tab:RmRcR}
\end{center}
\end{table}

\section{The Simplicial Ricci Flow Equation}
\label{sec:theorem} 
\noindent
{\bf Definition 1.}
We define the dual-edge Regge-Ricci flow equation for any compact, piecewise--flat simplicial geometry, $\mathcal S$, as an equation for each edge, $\lambda$, in the circumcentric dual lattice, ${\mathcal S}^*$,
\begin{equation}
\left( \begin{array}{c}Fractional\ rate\ of\ change\\ of\  a\  circumcentric\\ dual\ edge,\\  \frac{1}{\lambda} \frac{\partial \lambda}{\partial t}\\   \end{array}  \right) = 
-\left( \begin{array}{c}Ricci\ tensor\ associated\\ to\ the\ circumcentric\\ dual\ edge,\\ Rc_\lambda    \end{array}  \right).
\end{equation}
\begin{remark}
We have shown here and in \cite{AMM:2011} that the Ricci tensor is naturally associated with a dual edge, $\lambda\in{\mathcal S}^*$.  We consider this edge, $\lambda$, and the associated hybrid cell, $V_{\lambda\lambda^*}$. We also use the $\ell\lambda$--hybrid orthogonal basis introduced in Sec.~\ref{subsec:Rm}.  Additionally, Sec.~\ref{subsec:Rc} and \cite{AMM:2011} showed that the mixed Ricci tensor is diagonal  in $\lambda$. Motivated by Hamilton's  RF equation, Eq.~\ref{Eq:frcc}, with mixed indices in an orthogonal basis, then the  RRF equation associated with this dual edge is, 
\begin{equation}
\frac{{\dot g}_{\lambda\lambda}}{g_{\lambda\lambda}} = -2\, R^\lambda{}_\lambda = -2\, Rc_\lambda.
\end{equation} 
Using the $g_{\lambda\lambda}$ component in the covariant  metric for the $\ell\lambda$--hybrid basis, Eq.~\ref{Eq:metricdown}, we obtain the dual-edge RRF equation for $\lambda$,
\begin{equation}
\label{Eq:dsrf}
\frac{\dot \lambda}{\lambda} = - Rc_\lambda.
\end{equation}
\end{remark}

Since the geometry of $\mathcal S$ and ${\mathcal S}^*$ is completely determined by the edges, $\left\{ \ell_i \right\}\in {\mathcal S}$, and nothing more, and since there are ordinarily many more dual edges, $\lambda$, than there are simplicial edges, $\ell$, then the system of equations given by the dual-edge RRF equation, Eq.~\ref{Eq:dsrf}, will ordinarily be overdetermined and have no solution.  In other words each $\lambda$ is a function of the $\ell$'s. Therefore, one solution would be for us to project, or trace, these dual-edge RRF equation onto the simplicial edges. This can be done using an appropriate weighted average.

\noindent
{\bf Corollary 1.}
The simplicial Regge-Ricci flow equation for any compact, piecewise--flat simplicial geometry, $\mathcal S$, can be expressed as an equation for each edge, $\ell \in {\mathcal S}$,\begin{equation}
\left( \begin{array}{c}The\  volume\ averaged\\  fractional\ rate\ of\ change\\ of\  the\ dual\ edge,\\  \left< \frac{1}{\lambda} \frac{\partial \lambda}{\partial t}\right>_{\ell}   \end{array}  \right) = 
-\left( \begin{array}{c}Ricci\ tensor\ associated\ to\\ the\ simplicial\ edge,\\ Rc_\ell    \end{array}  \right).
\end{equation} 
\begin{proof}
We express the RRF equation by mapping the dual-edge RRF equation (Eq.~\ref{Eq:dsrf}) onto an edge $\ell$,
\begin{equation}
\label{Eq:cp2}
\langle \dot \lambda/\lambda \rangle_\ell = -\langle  Rc_\lambda \rangle_\ell.
\end{equation}
It suffices then to consider the right--hand side of this equation. Using the trace in \cite{McDonald:2012,M:1997,MM:2008} on $Rc_\lambda$, we note that 
\begin{equation}
\label{Eq:cp}
Rc_\ell V_\ell = \sum_{\lambda_{|\ell^*}} Rc_\lambda V_{\ell\lambda}.
\end{equation}
We also note that the $\ell-\lambda$ reduced hybrid volume is symmetric in its indices,
\begin{equation}
V_{\ell\lambda} = V_{\lambda\ell},
\end{equation}
and we can expand the $\ell$ hybrid volume into a sum of these reduced hybrid cells,
\begin{equation}
V_\ell = \sum_{\lambda_{|\ell^*}} V_{\ell\lambda}.
\end{equation}
Therefore, we see by the definition of the weighted average in Eq.~\ref{Eq:cp} that the right-hand side of Eq.~\ref{Eq:cp2} is simply 
\begin{equation}
\langle  Rc_\lambda \rangle_\ell = Rc_\ell,
\end{equation}
and is given explicitly by Eq.~\ref{Eq:Rcell}.
Applying the edge-based weighted average to the LHS of the equation as we just did for the RHS of the equation yields the simplicial RRF equation,
\begin{equation}
\label{Eq:srfelleqn}
 \left< \frac{1}{\lambda} \frac{\partial \lambda}{\partial t}\right>_{\ell} = - Rc_\ell.
\end{equation}
\end{proof}

\section{The RRF equation in 3-Dimensions: Two Simple Analytic Applications}
\label{sec:RRF3}

\begin{figure}
\centering
\includegraphics[height=3.0in]{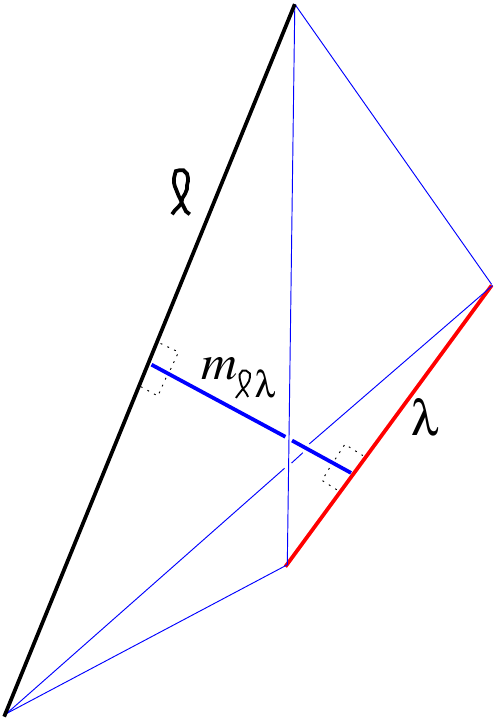}
\caption{We illustrate here the restricted hybrid volume, $V_{{\ell\ell^*}_{|\lambda}}$, which is just the usual hybrid volume $V_{\ell\ell^*}$ further restricted to one of the dual edges $\lambda$.  $\lambda$ is one of the dual edges in the boundary of dual polyhedron $\ell^*$. } \label{Fig:elllambda} 
\end{figure}

In 3-dimensions the RRF equation obtains a relatively simple form since the edges, $\ell$, are the hinges, $h$. Consider an edge $\ell_i\in S$ and the associated RRF equation, Eq.~\ref{Eq:srfelleqn}.  The right--hand side of this equation can be expressed in terms of the deficit angles of $\ell_j$, dual areas, $\ell^*_j$,  and certain ratios of reduced hybrid to hybrid volumes, 
\begin{equation}
RHS_i = -Rc_{\ell_i} = -2\, \sum_{{\lambda_k}_{|\ell^*_i}} \sum_{{\ell_j}_{|\lambda_k}}  \frac{\epsilon_{\ell_j}}{\ell^*_j}\, 
\left( \frac{V_{\lambda_k \ell_j}}{V_{\lambda_k}}\right) \left(\frac{V_{\ell_i \lambda_k}}{V_{\ell_i}} \right) .
\end{equation}
The left--hand side is a weighted average of the fractional rate of change of the dual edges, $\lambda_k \subset \ell^*_i$,
\begin{equation}
LHS_i =\sum_{{\lambda_k}_{|\ell_i^*}} \frac{\dot \lambda_k}{\lambda_k}\, \left(\frac{V_{\ell_i \lambda_k}}{V_{\ell_i}}\right),
\end{equation}
where $V_{\ell\lambda}$ is the reduced hybrid volume shown in Fig.~\ref{Fig:elllambda}, and $V_\ell = \frac{1}{3} \ell\ell^*$ is the volume of the hybrid block.  The mutual orthogonality of $m_{\ell \lambda}$, $\ell$, and $\lambda$ yields a relatively simple formula for the volume of this reduced hybrid tetrahedron, 
\begin{equation}
V_{\ell\lambda} = \frac{1}{6}\, \ell \lambda m_{\ell \lambda}.
\end{equation}
Therefore, the 3-dimensional RRF equation associated to edge, $\ell_j$, is simply 
\begin{equation}
\label{Eq:3srfsimple} 
\sum_{{\lambda_k}_{|\ell^*_i}} \, 
\left[
\frac{\dot \lambda_k}{\lambda_k} +
2 \sum_{{\ell_j}_{|\lambda_k}}  \frac{\epsilon_{\ell_j}}{\ell^*_j}\, 
\left( \frac{V_{\lambda_k \ell_j}}{V_{\lambda_k}}\right)
 \right]
\left(\frac{V_{\ell_i \lambda_k}}{V_{\ell_i}} \right)=0.
\end{equation}
Each dual edge, $\lambda_k\in \ell^*_i$,  is a function of 9 simplicial edges of the two tetrahedra in $\mathcal S$ that share the triangle $\lambda^*_k$; consequently , 
\begin{equation}
\dot \lambda_k = \sum_{{\ell_j}_{|\lambda_k^*}} \frac{\partial \lambda_k}{\partial \ell_j}\, \dot \ell_j.
\end{equation}
This yields a coupled set of first--order nonlinear algebraic equations, 
\begin{equation}\label{Eq:doublesumRRF3D}
\sum_{{\lambda_k}_{|\ell^*_i}} \sum_{{\ell_j}_{|\lambda_k}}
\left[
\frac{\partial \lambda_k}{\partial \ell_j}\, \frac{\dot \ell_j}{\lambda_k} + 2 \frac{\epsilon_{\ell_j}}{\ell^*_j}\, 
\left( \frac{V_{\lambda_k \ell_j}}{V_{\lambda_k}}\right)\right] \, \left(\frac{V_{\ell_i \lambda_k}}{V_{\ell_i}} \right)
= 0,
\end{equation}
 that should approximate the Hamilton RF dynamics in terms of the  evolution of the edges of the simplicial geometry.  In the next two subsections we will examine the solution of Eq.~\ref{Eq:3srfsimple} for two simple models, the 3--sphere and the 3--cylinder,  and demonstrate that they agree with their continuum analogues.

\subsection{5--Cell, 16--Cell and 600--Cell Model of 3--D RRF} 
\label{Sec:600}

In this section we will examine the RF of a 3--sphere geometry by comparing the continuum evolution with the solutions of the RRF equations (Eq.~\ref{Eq:3srfsimple}).  

The metric describing the geometry of a 3--sphere in spherical coordinates is 
\begin{equation}
ds^2 =  g_{ij}\, dx^i dx^j = a^2 \left( d\chi^2 + \sin\!^2{(\chi)} d\Omega^2\right);
\end{equation}
here, we let the radius of the sphere, $a(t)$,  depend on an external time parameter, $t$, and $d\Omega^2= d^2\theta+\sin{(\theta)} d\phi^2$ is the usual spherical line element. 
The only non--zero components of the mixed--component Ricci tensor for this geometry are
\begin{equation}
R^\chi\phantom{}_\chi = R^\theta\phantom{}_\theta = R^\phi\phantom{}_\phi =   (d-1)\, K = 2\, \frac{1}{a^2}.
\end{equation}
Therefore, Eq.~\ref{Eq:frcc} governing the RF for this geometry is a single equation for the radius, 
\begin{equation}
\label{Eq:cs}
\frac{\partial a^2(t)}{\partial t} = -4.
\end{equation}
Under RF this sphere collapses parabolically in time starting from an initial radius, $a_o = a(0)$ at time $t=0$,
\begin{equation}
\label{Eq:cpt}
a^2(t) = a^2_o - 4t.
\end{equation}

\begin{figure}
\centering
\includegraphics[height=1.6in]{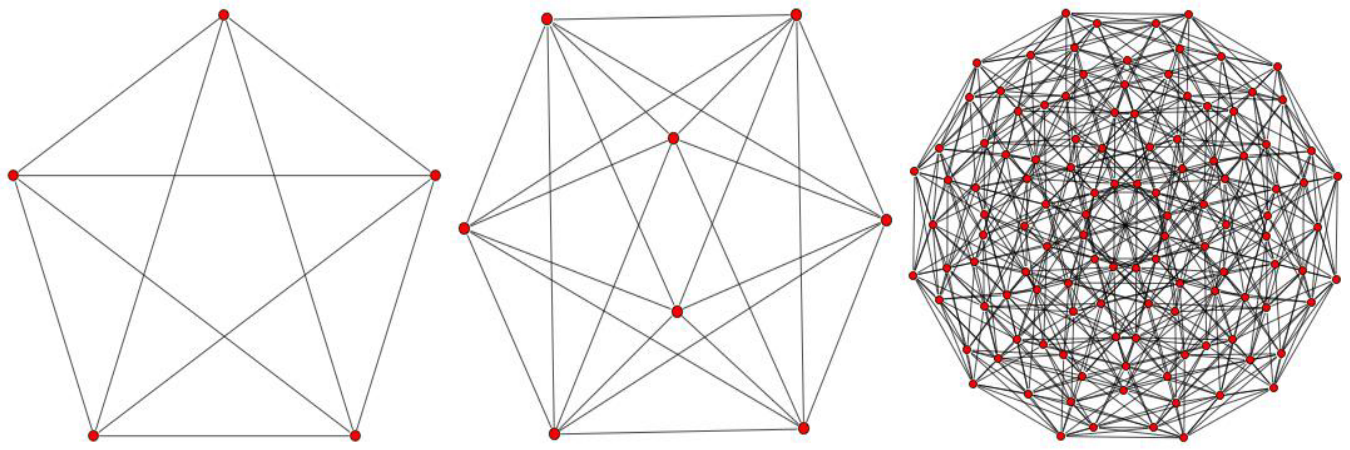}
\caption{We illustrate here a Coxeter plane representation for each of  the three regular 4--polytopes whose boundary are a simplicial approximation to the 3--sphere geometry.  The left 5--cell polytope is a simplex whose boundary consists of 5 regular tetrahedra, three of each sharing each of its 10 edges. The middle figure is an illustration of the 16--cell cross polytope.  It has 16 regular tetrahedra in its boundary. Each of its 24 edges is the meeting place of 4 of these regular tetrahedrons. Finally the right-most figure shows a Coxeter plane representation of the 600--cell model.  As its name suggests it has 600 regular tetrahedra in its 3--dimensional boundary.  Each of its 720 edges share 5 tetrahedra and in this sense this is often viewed as the higher-dimensional analogue of the icosahedron. This is the highest resolution homogeneous rendering of a 3--sphere geometry by a regular tetrahedral lattice.}\label{Fig:516600}
\end{figure}

We now analyze the RRF equations for three simplicial representations of $S^3$.  There are three regular polytopes in 4--dimensions \cite{Coxeter:1973} whose boundaries are built of regular tetrahedra and of which each is a homogeneous approximation to a 3--sphere.  These polytopes are referred  to as the simplex (5--cell model), the cross polytope (16--cell model) and the icosahedral pentagonal polytope (600--cell model), as depicted in Fig.~\ref{Fig:516600}.  The geometry of each of these three simplicial lattices are determined by a single length parameter, the length of the regular tetrahedron edge.  We could let the edges of these lattice geometries vary independently; however, for the purpose of this paper we assume that they are all equal and, therefore, the polytopes remain regular during their evolution under RF.

Here we explicitly construct  and solve Eq.~\ref{Eq:doublesumRRF3D} for each of these three models and compare their evolution with the continuum. In these models the curvature is concentrated at the edges. The 5--cell has three regular tetrahedra ($p=3$) sharing each edge, $\ell$,  so that the dual Voronoi polyhedron, $\ell^*$ is a triangle. The 16--cell model has four regular tetrahedra ($p=4$) sharing each edge, $\ell$,  so that the dual Voronoi polyhedron, $\ell^*$, is a square.  Finally, the 600--cell model has five regular tetrahedra ($p=5$) sharing each edge, $\ell$,  so that the dual Voronoi polyhedron, $\ell^*$ is a pentagon.  
\begin{equation}
\ell^*_p = \frac{p}{2}\, \lambda_p \, m_{{}_{\ell_p \lambda_p}} = \frac{\sqrt{2}\, p}{24}\, \ell^2_p.
\end{equation}
The dihedral angle of a regular tetrahedron is given by 
\begin{equation}
\theta = \sec^{-1}\left(3\right) \approx 1.23096,
\end{equation}
so that the deficit angle $\epsilon_p$ at edge $\ell_p$  for each of these models is simply, 
\begin{equation}
\epsilon_p = 2\pi - p\, \theta,
\end{equation}
where $p=3$, $4$ or $5$ for the 5--cell, 16--cell and 600--cell model, respectively. In particular, 
\begin{equation}
\begin{array}{llr}
5-\text{cell},\  & p=3,\ \  \epsilon_3 & \approx 2.59031 \\
16-\text{cell} ,\ & p=4,\ \   \epsilon_{4} & \approx 1.35935 \\
600-\text{cell},\  & p=5,\ \    \epsilon_{5} & \approx 0.12839.
\end{array}
\end{equation}
The dual edges, $\lambda_p$, for each of these models are proportional to the edge length,
\begin{equation}
\lambda_p = \frac{\sqrt{6}}{6}\, \ell_p.
\end{equation}
The moment arm, $m_p := m_{_{\ell_p \lambda_p}}$, in Eq.~\ref{Eq:3srfsimple} is defined with respect to an edge, $\ell_p$, and a dual edge, $\lambda_p$, in the Voronoi polygon, $\ell^*$.  This line segment reaches from the center of the edge, $\ell_p$,  to the circumcenter of the triangle, $\lambda^*_p$.  The moment arm, for each of our three lattices, is also proportional to the simplicial edge length,  
\begin{equation}
m_p = \frac{\sqrt{3}}{6}\, \ell_p.
\end{equation}
As all the dual edges, $\lambda_p$, are of equal length, and as all the simplicial edges, $\ell_p$ are of equal length in this model, then the dual-edge RRF equations are equivalent to the simplicial RRF equations. We therefor can consider here the left-hand side of the dual-edge RRF equation, 
\begin{equation}
LHS = \frac{\dot \lambda_p}{\lambda_p} = \frac{\dot \ell_p}{\ell_p},
\end{equation}
and its right-hand side,
\begin{equation}
RHS = -Rc_{\lambda_p} = -2\, \frac{\epsilon_p}{\ell^*_p}.
\end{equation}
These yield the dual-edge RRF equation for each of these three models, 
\begin{equation}
\label{Eq:peqn}
 \frac{\partial \left(\ell^2_p\right)}{\partial t}  = - \frac{48 \sqrt{2}}{p}\, \epsilon_p
\end{equation}

How do we compare this with the continuum; in other words, how do we relate the radius, $a$, with the simplicial edge, $\ell_p$? One way we  can supply such a relation is by equating the surface 3--volume of the boundary of each polytope with the 3--volume of the 3--sphere of radius, $a$; i.e. 
\begin{equation}
\underbrace{
\left(
\begin{array}{c}
\hbox{Volume of a 3--sphere}\\
\hbox{of radius}\ a
\end{array}
\right)}_{2\pi^2 a^3} 
= 
\underbrace{
\left(
\begin{array}{c}
\hbox{3--Volume of the boundary of}\\
\hbox{the $p$--cell with edge $\ell_p$}
\end{array}
\right)}_{ N_p\,\left(\frac{\sqrt{2}}{12}\ell^3_p\right)}.
\end{equation} 
Here, $N_p=5$, $16$ and $600$ is the number of tetrahedra in the boundary of the 5--cell, 16--cell and 600--cell polytope; respectively.
Therefore, we obtain the desired relation between the radius and the simplicial edge length, 
\begin{equation}
\label{Eq:ella}
\ell_p^2 = \left( \frac{24 \pi^2}{ \sqrt{2}\, N_p}  \right)^{2/3} a^2.
\end{equation}

So that we may compare the RRF evolution with the continuum solution (Eq.~\ref{Eq:cs}), we substitute Eq.~\ref{Eq:ella} into the RRF equation (Eq.~\ref{Eq:peqn}), yielding
\begin{equation}
\frac{\partial}{\partial t} \, a^2(t)  = -2\, \underbrace{\left\{\frac{6 \sqrt{2}}{p}\,\left( \frac{\sqrt{2}\, N_p}{3\, \pi^2} \right)^{2/3}\epsilon_p\right\}}_{\approx R_{\chi\chi} = 2}.
\end{equation}
We observe from this RRF equation that, as in the continuum  Eq.~\ref{Eq:cpt},  each of the three regular polytopes collapses parabolically in time starting from an initial effective radius, $a_o = a(0)$ at time $t=0$,
\begin{equation}
a^2(t) = a^2_o - 2\, \underbrace{\left\{\frac{6 \sqrt{2}}{p}\,\left( \frac{\sqrt{2}\, N_p}{3\, \pi^2} \right)^{2/3}\epsilon_p\right\}}_{\approx 2} t.
\end{equation}
In Table~\ref{Table:s3converge}, we show the percent deviation of the component of the Ricci tensor, $R_{\chi\chi}$,  from its continuum value $R_{\chi\chi}=2$ as a function of the deficit angle for each of these three models.  These results are consistent with the predicted second--order convergence of the simplicial Ricci tensor to its continuum value as a function of deficit angle.  In particular, the logarithm of this fractional deviation for the 16 to the 600--cell model  gives  a convergence rate of approximately second order,  ($slope \approx 1.88$).
 \begin{table}[h]
 \begin{center}
\begin{tabular}{| c || c | c | }
\hline
Model & Deficit Angle & Percent Deviation from the \\
& $\epsilon_\ell$  ( radians)  & Continuum Value of $R_{\chi\chi}=2$ \\
\hline\hline
5-cell & 2.59031 & 41.0\% \\
\hline
16-cell & 1.35935 & 20.5\% \\
\hline
600-cell & 0.12839 & 2.02\% \\
\hline
\end{tabular}
\caption{This table demonstrates the convergence of the RRF equations to their continuum counterpart as a function of  increasing mesh refinement.  The 600--cell model reproduces the collapse rate in the square of the radius of the sphere within $\sim2\%$.  The results from the 16--cell to the 600--cell model show a convergence rate in the deficit angle of order $\sim$2 which is consistent with the non--commutativity of the finite rotations,($\epsilon_\ell$), as was argued in \cite{M:1986}.  }\label{Table:s3converge}
\end{center}
\end{table}

We are currently examining the solutions of the RRF equations on a 600--cell model where we allow all of its 720 edges to freely vary\cite{Ray:2012}.  The RRF equations in this case form a sparsely--coupled set of 720 nonlinear first order differential equations for the 720 simplicial edge lengths. Our goal with this expanded model is threefold.  First, we wish to reproduce the exponential uniformization proven by Hamilton for those perturbations of the edge lengths that keep the curvature positive. Second, we wish to analyze the behavior and stability of the RRF equations derived here. Thirdly, we can begin examining the role that modified RF and simplicial diffeomorphism freedom have on the stability and the accuracy of RRF.  This model also has 1200 dual edges, $\lambda$, so that we can examine more closely the behavior of the projection of the dual RRF equations onto the simplicial edges, $\ell$. 

\subsection{Icosahedral 3--Cylinder Model of 3--D RRF}

In this section we apply our RRF equations to a 3--dimensional cylinder of radius, $R(t)$.  The metric for a 3--cylinder is given by, 
\begin{equation}
ds^2 = g_{ij} dx^idx^j = dz^2 + R^2\left(t\right)\, d\Omega^2 ,
\end{equation}
where $z\in\{-\infty,\infty\}$ is the axial coordinate, and $d\Omega^2 =  d\theta^2 + \sin{\left(\theta\right)} d\phi^2$ is the usual spherical line element. We let the radius, $R(t)$,  of the 2--sphere cross sections of the 2--cylinder depend on an external time parameter, $t$. For this metric, the only non--zero components of the mixed--component Ricci tensor for this geometry are,
\begin{equation}
R^\theta\phantom{}_\theta = R^\phi\phantom{}_\phi =   \frac{1}{R^2(t)}.
\end{equation}
Therefore, Eq.~\ref{Eq:frcc} governing the RF for this geometry is a single equation for the radius, 
\begin{equation}
\label{Eq:cs}
\frac{\partial R^2(t)}{\partial t} = -2.
\end{equation}
Under RF this sphere collapses parabolically in time starting from an initial radius, $R_o = a(0)$ at time $t=0$,
\begin{equation}
R^2(t) = R^2_o - 2t.
\end{equation}

\begin{figure}[h]
\centering
\includegraphics[height=1.25in]{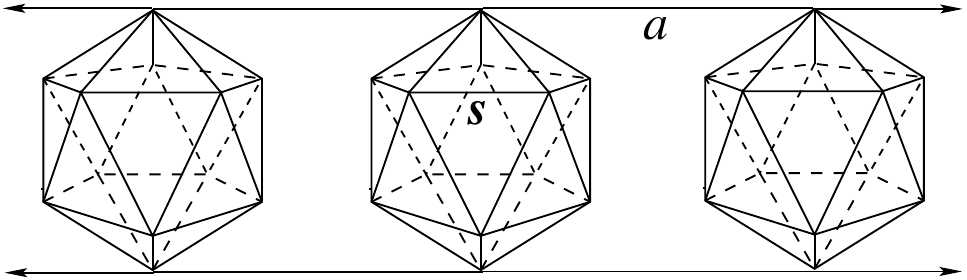}
\caption{This illustrates the icosahedral 3--cylinder model we use to approximate the RF of a 3--cylinder.  In this lattice, the spherical cross section of the 3--cylinder is approximated by a regular icosahedron with edge length, $s$.  Adjacent icosahedra are joined together at their corresponding vertices by 20 equal axial edges, $a$.  In this model we examine the evolution under the RRF equations of these two edge lengths. }\label{Fig:3cylinder}
\end{figure}

\begin{figure} [h]
\centering
\includegraphics[height=2.5in]{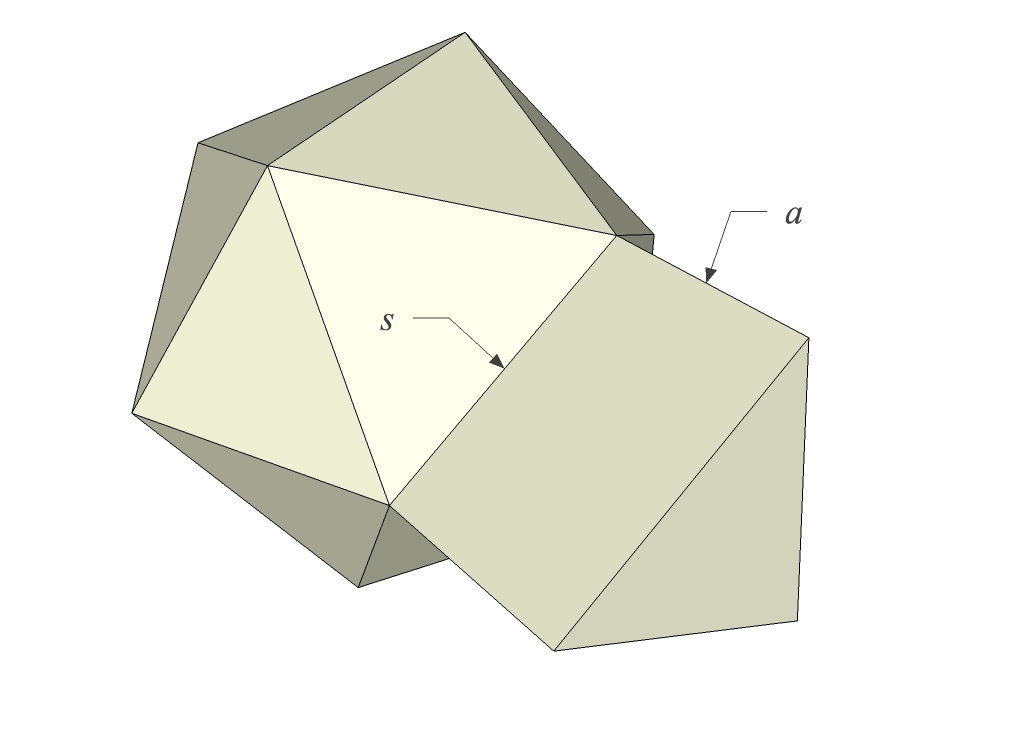}
\caption{We show here a 3--dimensional illustration of the cross--sectional icosahedron and, on one of its 20 triangular faces,  a regular triangular prism. This prism reaches from the icosahedron and will be attached to the corresponding triangular face of the adjacent icosahedron.  Our  lattice representation of a 3--cylinder has only one polyhedral block, namely this regular triangular prism. }\label{Fig:3Dcylinder}
\end{figure}

\begin{figure}[h]
\centering
\includegraphics[height=2.5 in]{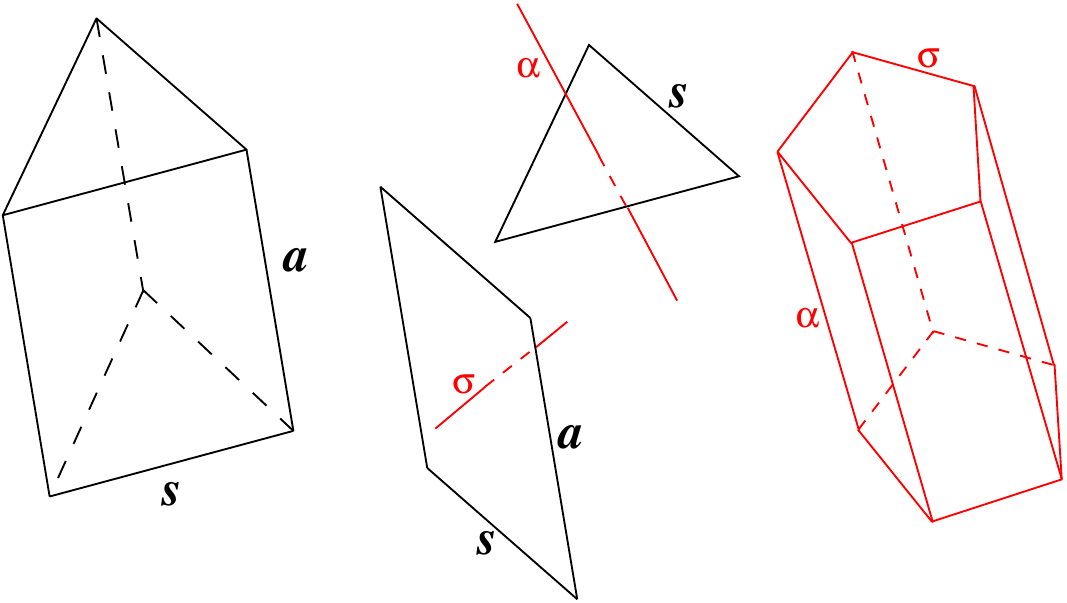}
\caption{We show on the left the triangular prism used to generate $\mathcal S$.  Due to the symmetry, its geometry is determined by the two edge lengths, $s$ and $a$. On the right we show the polyhedron lattice block in the circumcentric dual lattice, $\mathcal S^*$. Its geometry is also fixed by assigning edge lengths to $\sigma$ and $\alpha$. In the center we show the two polygons in $\mathcal S$ dual to the edges, $\alpha$ and $\sigma$  in $\mathcal S^*$.   }\label{Fig:PrismMoments}
\end{figure}

We now analyze the RRF equations, Eq.~\ref{Eq:doublesumRRF3D},  for a particular lattice representation of the 3--cylinder.  Similar to the last problem, the dual-edge RRF equations are equivalent to the simplicial RRF equations.   In particular we represent each of the 2--dimensional spherical cross sections of the 3--cylinder by an icosahedron with edge length $s(t)$.  Adjacent icosahedra are connected by axial edges of length $a(t)$ as illustrated in Fig.~\ref{Fig:3cylinder}.   Since adjacent icosahedra are joined together by 12 identical axial edges, $a$, then this 3--cylinder model consists of only one type of block -- a regular triangular prism shown in Fig.~\ref{Fig:3Dcylinder} and Fig.~\ref{Fig:PrismMoments}. That this block is not a tetrahedron is not problematic.  Ordinarily, the geometry of a non--simplicial block is not uniquely determine by its edge lengths, (i.e. it is non--rigid).  However, the geometry of this prism is determined uniquely by its  two  edge lengths, $s$ and $a$, because of its symmetry.  
  
The circumcentric dual lattice for this model is similar to a regular triangular prism lattice.  The cross sections are regular dodecahedra of edge length 
\begin{equation}
\sigma=(\sqrt{3}/3) s,  
\end{equation}
and adjacent dodecahedra are connected to each other by axial edges 
\begin{equation}
\alpha=a,
\end{equation}
as shown in Fig.~\ref{Fig:PrismMoments}.
The dual lattice is composed of a single block -- a regular pentagonal prism.   The  dodecahedrons are  located halfway between each icosahedron, and therefore the adjacent dodecahedrons are connected by 20 equal axial edges of length $\alpha=a$.  Each $\alpha$--edge is dual to a face of an icosahedron, while each $\sigma$ is dual to a rectangular face of the triangular prism.  

In this model, the only non--zero curvature is concentrated at the edges $a$.  As we will see, there is no curvature along edge $s$ .  In order to calculate the curvature at edge $a$ and edge $s$ we need to calculate each of  their deficit angles. We note that there are two dihedral angles of a regular triangular prism.  The dihedral angle associated with the triangle edge, $s$,  is a right angle, 
\begin{equation} 
\theta_s = \frac{\pi}{2},
\end{equation} 
while the dihedral angle associated to the axial edge, $a$, is equal to the internal angle of the regular triangle,
\begin{equation}
\theta_a = \frac{\pi}{3}.
\end{equation}
Each axial edge $a$ has five regular triangular prisms hinging on it, so the deficit angles for edges $a$ are 
\begin{equation}
\epsilon_a = 2\pi - 5 \theta_a = \frac{\pi}{3},
\end{equation} 
and since there are four prisms sharing each edge icosahedral edge $s$, the deficit angle vanishes,  
\begin{equation}
\epsilon_s = 2\pi - 4 \theta_s = 0.
\end{equation}
The only non-zero curvature is along the axial edges $a$.  

The two edges, $s$ and $a$,  and their respective dual rectangle $s^*$ and dual pentagon $a^*$,  define three distinct  moment arms. There are two moment arms associated with icosahedral edge, $s$.  One of these connects the midpoint of $s$ with the midpoint of the dual edge $\sigma$,
\begin{equation}
\label{Eq:ms}
m_{s\sigma} = \frac{a}{2},
\end{equation} 
and the other connects the midpoint of $s$ with the midpoint of the dual edge $\alpha$,
\begin{equation}
\label{Eq:mb} 
m_{s\alpha} = \frac{\lambda_s}{2}= \frac{\sqrt{3}}{6} s.
\end{equation}
There are two of these moment arms on the rectangle $s^*$.    
Finally, there are five identical moment arms associated with edge $a$.  These moment arms reach from the center of edge $a$ to the center of edge $\sigma$; hence, 
\begin{equation}
m_{a\sigma} = \frac{s}{2}.
\end{equation}

The left--hand side of the RRF equations, Eq.~\ref{Eq:doublesumRRF3D}, for this  model is given by equations, 
\begin{eqnarray}
LHS_a & = &  \sum_{\lambda_{|a^*}}  \frac{\dot \lambda}{\lambda} \, \left(\frac{V_{\ell \lambda}}{V_\ell}\right) =  \frac{ \dot \sigma}{\sigma} =\frac{\dot s}{s} , \\
LHS_s & = &  \sum_{\lambda_{|s^*}}   \frac{\dot \lambda}{\lambda} \, \left(\frac{V_{\ell \lambda}}{V_\ell}\right) =  \frac{\dot \sigma}{\sigma} +  \frac{\dot \alpha}{\alpha} = \frac{\dot s}{s} +  \frac{\dot a}{a},
\end{eqnarray}
and the right-hand side of these equations are simply, 
\begin{eqnarray}
RHS_a & = & -Rc_a = -20 \left(  \frac{\epsilon_s}{s^*} \left(  \frac{V_{\sigma s}}{V_\sigma} \right) + \frac{\epsilon_a}{a^*} \left( \frac{V_{\sigma a}}{V_\sigma}\right)  \right) \left( \frac{V_{a\sigma}}{V_a} \right) = -\frac{\epsilon_a}{a^*}\\
RHS_s & = & -Rc_s = -2\, \frac{\epsilon_s}{s^*} -  \frac{\epsilon_a}{a^2} = - \frac{\epsilon_a}{a^*} .
\end{eqnarray}

Consequently, this yields the two RRF equations for this icosahedral 3--cylinder.  
\begin{eqnarray}
\frac{\partial s^2}{\partial t}  & = &  -\frac{8\pi}{5\sqrt{3}} \label{Eq:srf3c-s}\\
\frac{\dot a}{a}  & = & 0. \label{Eq:srf3c-a}
\end{eqnarray}

How do we compare this with the continuum; in other words, how do we relate the radius, $r$ with the simplicial edge, $s$? One way we can supply such a relation is by equating the surface area of a  2--sphere of radius $r$ with the surface area of an icosahedron with edge length $s$; i.e. 
\begin{equation}
\underbrace{
\left(
\begin{array}{c}
\hbox{Area of a 2-sphere}\\
\hbox{of radius}\ r
\end{array}
\right)}_{4\pi r^2} 
= 
\underbrace{
\left(
\begin{array}{c}
\hbox{Surface area of an}\\
\hbox{icosahedron with edge $s$}
\end{array}
\right)}_{5\sqrt{3}\, s^2},
\end{equation} 
or, 
\begin{equation} 
s^2 = \left( \frac{4 \pi}{5 \sqrt{3}} \right)\,  r^2.
\end{equation} 
Substituting this into  Eq.~\ref{Eq:srf3c-s} gives us the RRF equation in terms of an effective radius, $r$, that is, serendipitously, the exact result as we found in the continuum, Eq.~\ref{Eq:cs},
\begin{equation}
\frac{\partial r^2}{\partial t} = - 2,
\end{equation}
and there is no change in time of the axial edge, Eq.~\ref{Eq:srf3c-a}, as is the case in the continuum . 

In 3--D and higher RF, singularities in the geometry may form.  The type--1 and type--2 neck-pinch singularities examined in the dumbbell and javelin models by \cite{Knopf:2004,Knopf:2011} can be examined in RRF by allowing the icosahedral edge lengths, $s$, and the axial edges, $a$,  to vary spatially  from one icosahedron to another.  The evolution of the dumbbell and javelin model is currently being examined \cite{Tison:2012} using the RRF equations introduced in this manuscript.  We have already observed  neck pinching qualitatively  in RRF, a more detailed comparison with the continuum models, analysis of conformal and diffeomorphic modes, and a stability analysis are currently underway \cite{Tison:2012}.    

\section{Future Directions: Singularities,  Surgery, Stability and the Simplicial Diffeomorphisms}
\label{Sec:Fini}

In this paper we have formulated a discrete form of the Hamilton RF equations for a piecewise flat simplicial geometry for any dimension, $d$.  We imagine that this simplicial formulation of RF will provide new insights both numerically and, more importantly, in representing the singularity and soliton structure on piecewise flat simplicial lattices.  It may be possible to use the RRF mathematics to explore singularity structure in three and higher dimensions in novel ways \cite{Cao:2011}.  Perhaps there are minimal simplicial structures characterizing RF singularities of various types, in other words we can ask, ``Given a RF singularity (Type-1, Type-2 etc.), what is the simplest simplicial lattice that exhibits this singularity?"

Hamilton demonstrated that the RF equations  will,  under certain conditions,  evolve to produce singularities  \cite{Hamilton:1995}.  He discussed some intuitive solutions for dumbbell geometries producing pinching singularities.  For the RRF equations to produce a faithful representation of Hamilton's RF, the lattice geometry will need to evolve so as to properly approach these pinching singularities and then to avoid them by surgery \cite{Perelman:2003,Hamilton:1995}.  The analysis of such singularities using the RRF equations is one of the primary directions guiding our current research.  We are  in the process of analyzing the dumbbell model with 3-sphere topology to reproduce the recent results of \cite{Knopf:2004,Knopf:2011} using the RRF equations.   Our preliminary results demonstrate neck pinching, although more analysis is necessary  \cite{Tison:2012}.  Once we demonstrate that the RRF equations reproduce quantitatively correct neck pinch singularity formation we will begin to examine a simplicial implementation of surgery.  This will involve separating the neck pinch into two separate lattices and then generate an appropriate simplicial cap for each of these, thereby continuing the RRF evolution to produce a connected sum of geometries. While we are only beginning in this direction, we believe it possible, at least in 3-dimensions,  to produce an automated surgery procedure for a lattice geometry in order to produce a complete RF evolution to a connected sum of Thurston-type geometries for a broad range of initial simplicial geometries and topologies \cite{Thurston:1997}.  

We have begun a numerical implementation of the RRF equations in 3-dimensions on 3-sphere lattices \cite{Ray:2012}.  In particular, we are currently analyzing the RRF of geometries with 3-sphere topology with everywhere positive Ricci curvature \cite{Ray:2012}.  In this analysis,  the initial simplicial geometry is a perturbation of the length of each of the 720 edges  of the 600-cell regular polytope (Sec.~\ref{Sec:600}) about their regular values so as to keep the $Rc$ positive everywhere.  Hamilton proved that non-uniformities in such geometries will converge to the regular 600-cell polytope exponentially fast \cite{Hamilton:1982} provided that the RRF faithfully reproduces continuum RF.  We are also examining Type-1 and Type -2 singularity formation in a simplicial representation of the dumbbell model \cite{Tison:2012,Gu:2012b}.  Our preliminary results for both models are qualitatively  consistent uniformization; however, we need to carefully examine the detailed exponential fast  uniformization predicted by Hamilton and to examine numerical behavior of the normalized RRF equations and their stability/instability.  For   small perturbation around the 600-cell model we examined the 720 eigenvalues of the Jacobian for the RRF equations and we find numerous positive real values indicating the well-known instability of the linearized RF equations about a 3-sphere geometry \cite{Ray:2012}.  Examining the stability of the RRF equations and providing suitable reformulations of these equations is the primary focus of our research.  

We are animated by three directions for this analysis.  First, Perelman introduced a gauge  fixed  (conformal and diffeomorphism) modified RF evolution  by introducing  an entropy integral -- an integral which produces uniformization via a  gradient flow \cite{Perelman:2003}.  We ask, ``What is the simplicial analogue of this for a piecewise flat simplicial lattice"? Second, an approximate simplicial diffeomorphism structure for a piecewise flat simplicial geometry has been understood \cite{GKMM:2009} and successfully implemented in RC.  This should enable us to examine a modified RRF equation using these approximate diffeomorphic degrees of freedom in order to improve the uniformization of the simplicial geometry (e.g. well-centered Delaunay lattices, maximally regular simplexes, etc.). Finally, there are many similarities between RF and issues involved in evolving black hole spacetimes that we believe can inform both fields of research. 

The numerical solution of Hamilton's RF equations in 3 and higher dimensions involve much the same issues as the vacuum black hole solutions of the Einstein equations.  In fact, the Einstein equations for a vacuum spacetime is just  the Ricci tensor.  First, Hamilton showed that one can encounter various singularities in RF which is also true for the numerical evolution of black hole spacetimes. Secondly, our recent results suggest that the RRF equations have instabilities in certain cases  \cite{Ray:2012,Tison:2012}  as do the discretized version of the standard ADM equations in numerical relativity.  Thirdly, Hamilton suggested a surgery procedure in order to integrate through the neck pinch singularities in much the same way as the moving puncture procedure or horizon excision procedure is used in numerical relativity to numerically evolve the spacetime geometry  avoid the black hole singularity.  Finally, in general relativity one utilizes the diffeomorphic degrees of freedom and conformal factors through the initial value data as well as a judicious choice of the lapse and shift functions during evolution so as to provide a stable solution of Einstein's equations  that also avoids singularities (e.g.  constant mean curvature slicing and minimal shear). In a similar vein, Perelman's entropy function utilizes diffeomorphic and conformal modes to provide a gradient flow for a modified RF equation.  This suggests that experience gained in general relativity may inform similar issues in the numerical analysis of RF, just as experience with RF evolutions may offer insights into numerical relativity.  We view this as an exciting avenue of research that may bridge two disciplines.

\section*{Acknowledgements:}
This work was supported in part under a grant from the Air Force Research Laboratory (AFRL/RITA), Grant \# FA8750-11-2-0089 and through support  from Air Force Office of Scientific Research through the American Society for Engineering Education's Summer Faculty Fellowship Program, and from Griffiss Institute under the AFRL VFRP.  JRM would like to acknowledge the support of the National Research Council.  We wish to thank Howard Blair, Matthew Corne,  Arkady Kheyfets, Shannon Ray, Chris Tison, Xuping Wang and Wei Zeng  for stimulating discussions and suggestions that improved this manuscript. Any opinions, findings and conclusions or recommendations expressed in this material are those of the author(s) and do not necessarily reflect the views of the AFRL.


\begin{thebibliography}{99}
\bibitem{Cao:2003}H-D. Cao, B. Chow, S-C Chu \& S-T Yau, eds., {\em Collected Papers on Ricci Flow}  in Series in Geometry and Topology, Volume 37 (International Press; Somerville, MA; 2003).
\bibitem{Chow:2004} B. Chow \& D. Knopf, {\em The Ricci Flow: An Introduction}, Mathematical Surveys and Monographs, Volume 110 (American Mathematical Society; Providence, RI; 2004).
\bibitem{Chow:2006} B Chow, P. Lu \& L. Ni, {\em Hamilton's Ricci Flow}, Graduate Studies in Mathematics, Volume 77 (American Mathematical Society; Providence, RI; 2006). 
\bibitem{Chow:2007} B. Chow, S-C Chu, D. Glickenstein, C. Guenther, J. Isenberg, T. Ivey, D. Knopf, P. Lu, F. Luo \& L. Ni, {\em The Ricci Flow: Techniques and Applications, Part 1: Geometric Aspects}, Mathematical Surveys and Monographs, Volume 135 (American Mathematical Society; Providence, RI; 2007).
\bibitem{Perelman:2003} G. Perelman, ``The entropy formula for the Ricci flow and its geometric applications,'' preprint,
math.DG/0211159; G. Perelman, ``Ricci flow with surgery on three-manifolds,'' preprint, math.DG/0303109; \& 
G. Perelman, ``Finite extinction time for the solutions to the Ricci flow on certain three-manifolds,''  preprint, math.DG/0307245.
\bibitem{Gu:2012} X. Yu, X. Yin, W. Han, J. Gao \& X. Gu, ``Scalable routing in 3D high genus sensor networks using graph embedding,'' {\em  INFOCOM 2012}: 2681-2685; 
Y. Wang, J.  Shi, X.  Yin, X.  Gu, T. F. Chan, S-T Yau, A. W. Toga, P. M. Thompson, ``Brain surface conformal parameterization with the Ricci flow,'' {\em IEEE Trans. Med. Imaging} {\bf 31}(2)  (2012) 251-264.
X. Gu, F.  Luo \& S-T Yau, ``Fundamentals of computational conformal geometry,'' {\em Mathematics in Computer Science} {\bf 4}(4) (2010) 389-429; 
B. Chow \& F. Luo, ``Combinatorial Ricci flows on surfaces,'' {\em J. Differential Geometry} {\bf 63} (2003) 97-129.
\bibitem{Peiro:2005} J. Peiro \& S. Sherwin, {\em Finite Difference, Finite Element and Finite Volume Methods For Partial Differential Equations,}  in Handbook of Materials Modeling, Volume 1,  Methods and Models, Springer, 2005.
\bibitem{Humphries:1997} S. Humphries, Jr., {\em Finite-Element Methods for Electromagnetism}, http://www.fieldp.com/freeware/finite\_element\_electromagnetic.pdf; originally published as {\em Field Solutions on Computers}  (ISBN 0-8493-1668-5) (Taylor and Francis, Boca Raton, 1997).
\bibitem{Regge:1961} T. Regge, ``General relativity without coordinates,'' {\em Il Nuovo Cimento} {\bf 19} (1961) 558-571.
\bibitem{Gentle:1998} A. P. Gentle, W. A. Miller, ``A fully (3+1)-D Regge calculus model of the Kasner cosmology,'' {\em Class.Quant.Grav.} {\bf 15} (1998) 389-405.
\bibitem{Hirani:2005} M. Desbrun, A. N. Hirani, M. Leok \&  J. E. Marsden, ``Discrete exterior calculus,''
e-print arXiv:math/0508341v2 [math.DG] on arxiv.org (2005).
\bibitem{Chow:2003} B. Chow \& F. Luo, ``Combinatorial Ricci flows on surfaces," {\em J. Diff. Geom.} {\bf 63} (2003) 97-129.
\bibitem{Glickenstein:2011} D. Glickenstein, ``Discrete conformal variations and scalar curvature on piecewise flat two- and three-dimensional manifolds," {\em J. Diff. Geom.} {\bf 87} (2011) 201-238.
\bibitem{G:2005} D. Glickenstein, ``Geometric triangulations and discrete Laplacians on manifolds,'' arXiv:math/0508188 [math.MG].
\bibitem{Ge:2013} H. Ge, ``Discrete Quasi-Einstein Metrics and Combinatorial Curvature Flows in 3-Dimension,'' 	arXiv:1301.3398 [math.DG].
\bibitem{Forman:2003}R. Forman,  ``Bochner's Method for Cell Complexes and Combinatorial Ricci Curvature,'' {\em Discrete Comput. Geom.} {\bf 29} (2003) 323Ð374.
\bibitem{LinYau:2010}Yong Lin and Shing-Tung Yau, ``Ricci curvature and eigenvalue estimate on locally finite graphs,'' {\em Math. Res. Lett.} {\bf 17} (2010) 343Ð356.
\bibitem{Knopf:2009}  D. Knopf, ``Estimating the trace-free Ricci tensor in Ricci flow,''  {\em Journal: Proc. Amer. Math. Soc.} {\bf 137} (2009), 3099-3103.
\bibitem{Thurston:1997} W.  Thurston, {\em Three-dimensional geometry and topology,} Vol. 1. Edited by Silvio Levy, Princeton Mathematical Series, 35,  (Princeton University Press, Princeton, NJ, 1997). 
\bibitem{McDonald:2012} J. R. McDonald, W. A. Miller, P. M Alsing, X. D. Gu, X. Wang \& S-T Yau, ``On exterior calculus and curvature in piecewise-flat manifolds,'' {\em paper submitted to J. Math. Phys.} (2012) arxiv.org/abs/1212.0919.
\bibitem{Hamilton:1982} R. Hamilton, ``Three-manifolds with positive Ricci curvature,'' {\em J. Diff. Geom} {\bf 17} (1982), 255-306.
\bibitem{Cartan:1951}  E. Cartan, Le\c ons sur la G\' eom\' etrie des Espaces de Riemann (Gauthier-Villars; Paris, 1951); translated into English by J. Glazebrook, Geometry of Riemann Spaces, Lie Groups: History Frontiers and Applications {\bf XIII}  (Math Science  Press, Brookline, MA, 1983).
\bibitem{FL:1984} R. Friedberg and T. D. Lee, {\em Nucl. Phys.} {\bf B242} (1984) 145-166.
\bibitem{M:1986} W. A. Miller, ``The geometric content of the Regge equations as illuminated by the boundary of a boundary principle," {\em Found. Phys.} {\bf 16} (1986) 143-169.
\bibitem{M:1997} W. A. Miller, ``The Hilbert action in Regge calculus,'' {\em Class. Quantum Grav.} {\bf 14} (1997) L199-L204. 
\bibitem{MM:2008} J. R. McDonald \& W. A. Miller, ``The scalar curvature in Regge calculus,'' {\em Class. Quantum Grav.} {\bf 25} (2008) 195017.
\bibitem{AMM:2011} P. M. Alsing, J. R. McDonald \& W. A. Miller, ``The simplicial Ricci tensor," {\em Class. Quantum Grav.} {\bf 28} (2011) 155007 (17 pp).
\bibitem{Van:2009} E. Vanderzee, A. N. Hirani, D. Guoy \& E. A. Ramos, arXiv:0802.2108v3 [cs.CG] 18 Aug 2009.
\bibitem{Caselle:1989}{M. Caselle, A. D'Adda \& L. Magena}, {\em Phys. Lett.} {\bf B232} (1989) 457.
\bibitem{MM:2008b} J. R. McDonald, W. A. Miller, ``A discrete representation of Einstein's geometric theory of gravitation: the fundamental role of dual tessellations in Regge calculus,''  submitted to "Tessellations in the Science: Virtues, Techniques and Applications of Geometric Tilings," ed. R. van de Weijgaert, G. Vegter, J. Ritzerveld and V. Icke (2008);  arXiv:0804.0279
\bibitem{Hartle:1981} J. B. Hartle \& R. Sorkin, ``Boundary terms in the action for the Regge calculus,''  {\em General Relativity and Gravitation}, {\bf 13}, no. 6 (1981) 541-549.
\bibitem{Sorkin:1975} R. Sorkin, "Time Evolution Problem in Regge Calculus." {\em Phys. Rev.} {\bf D12} (1975) 385-396.
\bibitem{Eisenhart:1997} L. P. Eisenhart, {\em Riemannian Geometry,} 8$^th$ edition (Princeton Univ. Press, Princeton, NY, 1997).
\bibitem{Lee:1984} R. Friedberg \& T.  D. Lee, ``Derivation of Regge's action from Einstein's theory of general relativity,''
 {\em Nuclear Physics}  {\bf B242} (1984) 145Ð166.
\bibitem{Coxeter:1973} H. S. M. Coxeter, {\bf Regular Polytopes}, (Dover Publ., NY, 1973).
\bibitem{Ray:2012} P. M. Alsing, M. Corne, J. R. McDonald, W. A. Miller, S. Ray, ``3D simplicial Ricci flow of the 600-cell model,'' {\em in preparation} (2013).
\bibitem{Tison:2012} P. M. Alsing, J. R. McDonald, W. A. Miller, C. Tison \& X. Wang, ``3D simplicial Ricci flow of the dumbbell model," {\em in preparation} (2013).
\bibitem{Knopf:2004} S.  Angenent \& D. Knopf, ``An example of neckpinching for Ricci flow on $S^{n+1}$,''  {\em Math. Res. Lett.} {\bf 11} (2004) 493-518.
\bibitem{Knopf:2011} S. Angenent, J. Isenberg and D. Knopf, ``Formal matched asymptotics for degenerate Ricci flow neckpinches,'' {\em Nonlinearity} {\bf 24} (2011), 2265-2280.
\bibitem{Cao:2011} H.-D. Cao, {\em Private communication,} UBC, Vancouver BC (2011).
\bibitem{Hamilton:1995} R. S. Hamilton, ``The formation of singularities in the Ricci flow,'' {\em Surveys in Differential Geometry} {\bf 2} (1995) 7-136.
\bibitem{Gu:2012b} X. Gu \& W. Zeng, {\it Private communication}, SUNY Stoney Brook (2012).
\bibitem{GKMM:2009} A. P. Gentle, A. Kheyfets, J. R. McDonald \& W. A. Miller, ``A Kirchhoff-like conservation law in Regge calculus,'' {\em Class. Quant. Grav.}  {\bf 26} (2009) 015005.
\end{thebibliography}
\end{document}